\documentclass{amsproc}
\usepackage{amssymb}
\usepackage{graphicx}
\usepackage{amscd}
\usepackage{amsmath}
\usepackage{graphicx}
\theoremstyle{plain}

\newtheorem{agreement}{Agreement}

\newtheorem{conclusion}{Conclusion}

\newtheorem{definition}{Definition}

\newtheorem{lemma}{Lemma}

\newtheorem{procedure}{Procedure}

\newtheorem{proposition}{Proposition}
\newtheorem{remark}{Remark}

\newtheorem{theorem}{Theorem}

\numberwithin{equation}{section}

\begin{document}

\title[Discrete dynamical systems]{Discrete dynamical
systems: inverse problems and related topics.}

\author{Alexander Mikhaylov }
\address{St. Petersburg   Department   of   V.A. Steklov    Institute   of   Mathematics
of   the   Russian   Academy   of   Sciences, 7, Fontanka, 191023
St. Petersburg, Russia and Saint Petersburg State University,
St.Petersburg State University, 7/9 Universitetskaya nab., St.
Petersburg, 199034 Russia.} \email{mikhaylov@pdmi.ras.ru}

\author{Victor Mikhaylov}
\address{St.Petersburg   Department   of   V.A.Steklov    Institute   of   Mathematics
of   the   Russian   Academy   of   Sciences, 7, Fontanka, 191023
St. Petersburg, Russia } \email{ftvsm78@gmail.com}

\date{July, 2019}

\begin{abstract}
In this review, we extend the Boundary Control method\,---\,an
approach to inverse problems based on control theory for dynamical
systems \,---\,to inverse problems for discrete dynamical systems.
We apply our results to classical moment problems, Toda lattices,
Weyl functions, de Branges spaces, Krein-Stieltjes strings, and
also to problems of numerical simulations.
\end{abstract}

\maketitle

\tableofcontents

\section{Introduction}

The theory of one-dimensional inverse dynamic problems for
hyperbolic systems is well-developed and has numerous applications
\cite{Blag2, KKL,Kab,R_VG}. However, the corresponding theory for
discrete systems is largely underexplored in the literature. In
this overview, we aim to summarize our contributions to inverse
problems (IP) involving one-dimensional discrete hyperbolic
systems, addressing both discrete and continuous time, as well as
related challenges in classical and spectral analysis.

The main mathematical tool we use in our research is the Boundary
Control (BC) method \cite{BDAN87,B97,B07}.  This method explores
connections between inverse problems in mathematical physics,
functional analysis, and control theory for partial differential
equations. It provides a powerful and reliable alternative to
traditional identification techniques that rely on spectral or
scattering methods. Initially developed to solve boundary inverse
problems in multidimensional wave equations, the BC method has
been successfully applied to various types of linear equations in
mathematical physics. For more detailed information, please refer
to the review papers \cite{B97,B07}, and the monograph \cite{KKL}.
Our work can be seen as an advancement and application of the BC
method to the field of discrete dynamical systems.

The main focus of our analysis is a Jacobi matrix. Consider a sequence of positive numbers
$\{a_0, a_1,\ldots\}$ and a sequence of real numbers $\{b_1, b_2,\ldots \}$. The Jacobi matrix
$A$ is then defined as follows:

\begin{equation}
\label{Jac_matr}
A=\begin{pmatrix} b_1 & a_1 & 0 & 0 & 0 &\ldots \\
a_1 & b_2 & a_2 & 0 & 0 &\ldots \\
0 & a_2 & b_3 & a_3 & 0 & \ldots \\
\ldots &\ldots  &\ldots &\ldots & \ldots &\ldots
\end{pmatrix}.
\end{equation}
The $N\times N$ block of this matrix is denoted by $A^N$.

Let $\mathbb{N}_0:=\mathbb{N}\cup \{0\}$. In the second section we
associate a discrete-time dynamical system with the matrix $A$ and
positive value $a_0$:
\begin{equation}
\label{Jacobi_dyn} \left\{
\begin{array}l
u_{n,\,t+1}+u_{n,\,t-1}-a_{n}u_{n+1,\,t}-a_{n-1}u_{n-1,\,t}-b_nu_{n,\,t}=0,\quad n,\in \mathbb{N},\, t\in \mathbb{N}_0\\
u_{n,\,-1}=u_{n,\,0}=0,\quad n\in \mathbb{N}, \\
u_{0,\,t}=f_t,\quad t\in \mathbb{N}_0.
\end{array}\right.
\end{equation}
We consider the complex sequence $f=(f_0,f_1,\ldots)$ as a
\emph{boundary control}. The solution to (\ref{Jacobi_dyn}) is
denoted by $u^f_{n,\,t}$. For a fixed  $\tau\in \mathbb{N}$, we
introduce the \emph{response operator} associated  with
(\ref{Jacobi_dyn}), which maps the control $f=(f_0,\ldots
f_{\tau-1})$ to $u^f_{1,\,t}$:
\begin{equation}
\label{resp_oper} \left(R^\tau f\right)_t:=u^f_{1,\,t},\quad
t=1,\ldots, \tau.
\end{equation}
The IP consists in recovering the sequences
$\{b_1,b_2,\ldots,b_n\}$, $\{a_0,a_1,\ldots,a_n\}$ for some $n\in
\mathbb{N}$ from $R^{2n}$. We also consider the case of a complex
matrix $A$. First section is central to the paper, containing key
formulas and operator representations that are utilized in the
subsequent sections. Detailed results of this section can be found
in \cite{MM1,MM4,MM18,MM19}. It is important to note that in
solving the dynamic inverse problem, we can work exclusively with
real-valued controls. However, the possibility to use complex
controls becomes significant in Sections
\ref{Section_Complex_matr} and \ref{Section_DeBranges}.

The classical moment problem consists in determination a Borel
measure $d\rho$ from the given a sequence of real numbers
$s_0,s_1,s_2,\ldots$, such that
\begin{equation}
\label{Moment_eq} s_k=\int\limits_{-\infty}^\infty
\lambda^k\,d\rho(\lambda),\quad k=0,1,2,\ldots
\end{equation}
holds. In the case a measure satisfying (\ref{Moment_eq}) exists,
these numbers referred to as \emph{moments}. It is well
established \cite{AH,Schm,S,T} that the spectral theory of Jacobi
operators is closely linked to moment problems. In the third
section, we employ inverse dynamic problems for discrete systems
as tools to tackle classical moment problems, specifically
focusing on questions of existence and uniqueness. Additionally,
we demonstrate that the solution to a moment problem can be
obtained as a weak limit of solutions to certain
finite-dimensional subproblems, which can be computed explicitly.
This section is based on the work presented in \cite{MM5,MM6}.

The semi-infinite Toda lattice can  be represented  \cite{To,Te} as the following infinite-dimensional nonlinear system:
\begin{equation}
\label{Toda_eq}
\begin{cases}
\dot a_n(t)=a_n(t)\left(b_{n+1}(t)-b_n(t)\right),\\
\dot b_n(t)=2\left(a_n^2(t)-a_{n-1}^2(t)\right),\quad t\geqslant
0,\, n=1,2,\ldots,
\end{cases}
\end{equation}
for which we seek for a solution that satisfies the initial
condition
\begin{equation}
\label{Toda_init} a_n(0)=a_n^0,\quad b_n(0)=b_n^0,\quad
n=1,2,\ldots.
\end{equation}
where $a_n^0,$ $b_n^0$ are real and $a_n^0>0,$ $n=1,2,\ldots$. The
methods for computing the functions $a_n(t)$, $b_n(t)$ have been
studied in detail, see for example \cite{Te,To} and the references
therein. In these studies, the authors employed techniques that
impose significant restrictions on the initial conditions, with
the most common assumption being that $a_n^0,$ $b_n^0$
$n=1,2,\ldots$ are bounded. However, the question of whether a
solution to (\ref{Toda_eq}) and (\ref{Toda_init}) can be
constructed for ``unbounded'' initial data remains important
\cite{IV1,IV2}.

In the fourth section, we present results from
\cite{MM12,MM13,MM17}, where the authors utilize the Moser formula
\cite{M75} to describe the evolution of the moments of the
spectral measure associated with the Jacobi matrix corresponding
to (\ref{Toda_eq}) and construct a solution to (\ref{Toda_eq}) and
(\ref{Toda_init}) for a wide class of unbounded initial sequences.

De Branges spaces play a crucial role in the inverse spectral
theory of  one-dimensional canonical systems \cite{DBr,DMcK,RR}.
In \cite{MM3}, the authors proposed a method to associate a de
Branges space with a specific dynamical system by considering the
initial-boundary value problem and employing the tools and
operators of the BC method. For the system described by
(\ref{Jacobi_dyn}), the algorithm outlined in \cite{MM3} is as
follows: by fixing a finite time  $t=T$ we introduce the
\emph{reachable set} of the dynamical system at this time:
\begin{equation*}
\mathcal{U}^T:=\{u^g_{\,\cdot\, ,T}\,|\, g\in \mathbb{C}^T\}.
\end{equation*}
Next, we apply the Fourier transform associated with the operator
corresponding to the matrix $A$ to elements from $U^T$, resulting
in a linear manifold  $\mathcal{F}U^T$.  This linear manifold is
then equipped with a norm defined by the \emph{connecting
operator}  $C^T$ associated with the system (\ref{Jacobi_dyn}),
resulting in a de Branges space. It's important to note that our
approach diverges from classical methods and may allow for
generalization to multidimensional systems \cite{MM23}.

Due to the finite speed of wave propagation in (\ref{Jacobi_dyn})
the procedure described above\,---\, precisely for dynamical
systems associated with Jacobi matrices\,---\, provides a
finite-dimensional de Branges space.  In contrast, the
construction of an infinite de Branges space corresponding to the
whole matrix $A$, rather than just the block $A^T$ for some $T\in
\mathbb{N}$ is addressed differently \cite{MM21}.

The fifth section provides outlining the construction of both
finite and infinite de Branges spaces for the system described in
(\ref{Jacobi_dyn}), following primarily \cite{MM3,MM21}. Readers
interested in additional examples of constructing de Branges
spaces for various dynamical systems using this method can refer
to \cite{MM15,MM16,MM21,MM23}.

In Sections \ref{Oper_fin} and \ref{Section_spectr_repr} we define
the operators corresponding to matrices $A^N$ and $A$. A key
object associated with these operators is the  \emph{Weyl
function} \cite{GS1,S1}. In the Section \ref{Moment_problem} we
utilize the dynamical system linked to $A^N$ to extract the
spectral data from dynamic data using the technique of the BC
method. The sixth section presents another significant
relationship between spectral and dynamic data:  by examining the
discrete dynamical systems described in  (\ref{Jacobi_dyn}) and
one corresponding to $A^N$  we derive a representation for the
Weyl function associated with the operators $A$ and $A^N$. For
more examples of such relationships, the reader can refer to
\cite{AMR,MM2,MM22}. The results presented in this section are
based on \cite{MMS}.

In the seventh section, we examine the continuous-time dynamical
system associated with $A^N$. Let $u=(u_1,\ldots,u_N)\in
\mathbb{R}^N$ and $T>0$ be fixed. We associate the following
dynamical system with the matrix $A^N$:
\begin{equation}
\label{DnSst}
\begin{cases}
u_{tt}(t)-A^Nu(t)=F(t),\quad t>0,\\
u(0)=u_t(0)=0,
\end{cases}
\end{equation}
where the vector function $F(t)=(f(t),0,\ldots,0)$, $f\in
L_2(0,T)$  is interpreted as a \emph{boundary control}. The
solution of (\ref{DnSst}) is denoted by $u^f$. We associate a
\emph{response operator}  with the system (\ref{DnSst}) which acts
according to the following rule:
\begin{equation}
\label{RespJM} \left(R^Tf\right)(t)=u^f_1(t),\quad 0<t<T.
\end{equation}
The dynamic IP for (\ref{DnSst}) consists of the reconstruction of
the matrix $A^N$ from the response operator $R^T$. It is important
to note that the properties of the systems (\ref{Jacobi_dyn}) and
(\ref{DnSst})  differ significantly regarding boundary
controllability and wave propagation speed, in (\ref{DnSst}) this
speed is infinite. These distinctions considerably  complicate the
solution of the IP because the standard tools of the BC method
\cite{B97,B07}, which rely on finiteness of wave propagation
speed, controllability, and geometrical optics, are not applicable
in this context. In this section, we outline the solution to the
IP and describe the characterization of dynamic inverse data. We
also explore the dynamic IP setting for a Krein string, a
historically significant object. We focus specifically on the case
of the finite Krein-Stieltjes string. Finally, we demonstrate how
a finite wave propagation speed emerges when approximating a
string with constant density using a weightless string with evenly
distributed masses. This section is based on
\cite{MM9,MM10,MM11,MM14}.

It is a well-known fact that Jacobi operators can be regarded as
discrete analogs of Sturm-Liouville operators \cite{Te}.
Consequently, it is natural to apply the concepts developed in
Section \ref{Section_Jacobi} to the treatment of both forward and
inverse numerical problems for ordinary differential equations
(ODEs), see also \cite{KabSh,KabSh1} for the application of some
special discretization in the BC method to the inverse acoustic
problem. In the final section we consider two such problems:
specifically, we explore applications in modeling forward and
control problems for the wave equation on quantum graphs, as well
as the numerical solution of a dynamic inverse problem for a
parabolic equation defined on a finite interval. We begin by
outlining these two problems in their continuous form and discuss
the rationale for transitioning to discrete models.

Let $ \Omega $ be a finite, connected, compact graph  consisting
of edges $E = \{e_1, \ldots, e_{N}\}$ that are connected at
vertices $V = \{v_1 \ldots, v_{M}\}$. Each edge $ e_j \in E $
corresponds to an interval $(0, l_j) $ of the real line. The
boundary $ \Gamma = \{v_1, \ldots, v_m \} $ of $ \Omega $
comprises the vertices with degree one (outer nodes). Denote by
$E(v)$ the set of edges incident to a vertex $v$. In what follows,
we assume that some boundary vertices are clamped, while a
non-homogeneous Dirichlet boundary condition is applied to the
remaining vertices.

We fix some $T>0$. The space of real square-integrable functions
on the graph $ \Omega $ is denoted by $L_2(\Omega ):=
    \bigoplus_{i=1}^{N}L_2(e_i).$ For the function $u \in L_2(\Omega)$, we write
    \begin{equation*}
        u:= \left\{ u^i\right\}_{i=1}^N\ ,\ u^i\in L_2(e_i).
    \end{equation*}

The forward problem is
\begin{align}
u^i_{tt}(x,t)-u^i_{xx}(x,t)=0,& \quad x \in e_i,\quad 0\leqslant t\leqslant T,\label{Wave_eqn} \\
u^i(v,\cdot)=u^j(v,\cdot),&  \quad e_i,\ e_j\in E(v),\quad v\in V\backslash\Gamma, \label{Cont}\\
\sum_{ e_i\in E(v)}\partial u^i(v,\cdot)=0,&\quad v\in V\backslash \Gamma, \label{Kirch}\\
u(v,t) = f(t),&\quad v\in\Gamma, \quad 0\leqslant t\leqslant T,\label{Bound}\\
u(x,0) = 0, u_t(x,0) = 0,& \ \quad  x\in\Omega. \label{Init}
\end{align}
Here, $\partial u^i(v,t)$ is the derivative of $u^i$ at point $v$
along the direction $e_i$, moving away from  $v$. This system
models small transverse oscillations of the graph, with $u(x,t)$
denoting the displacement of the point $x\in\Omega$  from its
equilibrium position at time $t$. It is assumed that at the
initial moment, the graph is stationary in an equilibrium state
and at the boundary of the graph the Dirichet boundary condition
is imposed. We note that when $f\in C^2(\Gamma)$, the solution to
(\ref{Wave_eqn})--(\ref{Init}) is classical.

In recent decades, forward and inverse problems on graphs have
gained significant attention. Studies on more general
systems\,---\,those with nonconstant density or
potential\,---\,have been conducted in references \cite{B_Tree,
BV, AK, AMN,Ku}. However, the algorithms proposed in these works
may be may be difficult to implement numerically. In particular,
after reconstructing the unknown coefficient on the edge, the
inverse data may be recalculated, which may lead to error
accumulation throughout the entire recursive process. To deal with
this the authors in \cite{AChMM} suggested focusing on a discrete
model instead of the continuous wave equation
(\ref{Wave_eqn})--(\ref{Init}), specifically the discrete wave
equation on a discrete quantum graph. In the first subsection, we
present the results of applying the variational principle to
derive a discretized version of the Kirchhoff conditions, which
most accurately reflect the wave propagation of the continuous
model. This work provides a discrete analog of equations
(\ref{Wave_eqn})--(\ref{Init}).

Let $T>0$ be a fixed parameter, and consider the interval $(0,L)$,
$L\in \mathbb{R}_+$. Given a potential $q\in C^2(0,L)$, we examine
the following initial-boundary value problem:
\begin{equation}
\label{Parab_eqn}
\begin{cases}
v_t(x,t)-v_{xx}(x,t)+q(x)v(x,t)=0,\quad 0<x<L, \, t\geqslant 0,\\
v(0,t)=f(t),\, v(L,t)=0,\quad  0<t<T,\\
v(x,0)=0,\quad 0<x<L.
\end{cases}
\end{equation}
Here, $f\in L_{2}(0,T)$ is a Dirichlet \emph{boundary control}. We
denote the solution to system (\ref{Parab_eqn}) as $v^f$. It is
known \cite{ABR} that $v^f\in C(0,T;H^{-1}(0,L))$.

For the system described by (\ref{Parab_eqn}), we introduce the
response operator $R^T:L_2(0,T)\mapsto L_2(0,T)$ with the domain
$D=\left\{f\in C^\infty(0,T)\,|\, f^{(k)}(0)=0,\,
k=0,1,\ldots\right\}$, defined by the
\begin{equation*}
\left(R^T\right)(t):=v^f_x(0,t),\quad t\in (0,T).
\end{equation*}
The dynamic inverse problem for system (\ref{Parab_eqn}) is to
determine the value of $q$ on the interval $(0, L)$ using the
known response operator $R^T$.

In the paper \cite{ABR}, the authors use the spectral
controllability of (\ref{Parab_eqn}) and the variational principal
to extract the spectral data of the operator
$-\frac{d^2}{dx^2}+q(x)$ on the interval $(0,L)$ with Dirichlet
boundary conditions from the dynamic data. After that, the IP can
be solved using standard methods.  While the time $T$ can be taken
arbitrary small in this approach, practical implementation seems
to be very challenging, and we do not know any attempts on the
numerical implementation of this method. To overcome this, we
replace the continuous system in (\ref{Parab_eqn}) with a discrete
version, resulting in a system of the form (\ref{Jacobi_dyn}). In
this section, we will discuss the results related to the solution
of the inverse problem and highlight the connections with the
results from Sections \ref{Section_Jacobi}, \ref{Oper_fin}, and
\ref{Section_tranc}. The results presented in this section are
based on the papers \cite{MM20} and \cite{MMSh}.

\section{Inverse dynamic problems for Jacobi matrices}
\label{Section_Jacobi}

In this section, we explain how to solve the dynamic inverse
problem for a dynamical system with discrete time, which is
associated with a Jacobi matrix $A$. First, we introduce the
initial-boundary problem and describe the setup for the inverse
problem. Next, we derive an important representation of the
solution and introduce the operators used in the BC method. Then,
we outline two methods for solving the inverse problem: one based
on Krein equations and another using factorization. We also
discuss the characterization of inverse data and consider a
special case involving a discrete Schr\"odinger operator. Finally,
we stdy the IP for a complex Jacobi matrix. Overall, this section
demonstrates how the BC method can be applied to both self-adjoint
and non-self-adjoint systems in discrete settings.

\subsection{Setup of the problem}

We associate a Jacobi matrix $A$ and a positive number $a_0$ with
a dynamical  system that has discrete time (equation
(\ref{Jacobi_dyn})). This system is a discrete version of the
initial-value problem for a wave equation with Dirichlet boundary
control at the endpoint on a semi-axis. When comparing this
discrete system to the continuous case, we notice that the
discrete analogue of the derivative with respect to $x$ at $x = 0$
is given by $u^f_{1,t} - u^f_{0,t}$. Therefore, the response
operator $R^{\tau}$ (equation (\ref{resp_oper})) is a discrete
analogue of a dynamic Dirichlet-to-Neumann operator.

The IP involves recovering sequences $\{b_1, b_2, \ldots, b_n\}$
and $\{a_0, a_1, \ldots, a_n\}$ for a given $n$ from the operator
$R^{2n}$. This problem serves as a discrete counterpart of an
inverse problem related to a wave equation on a semiaxis, where
the dynamic Dirichlet-to-Neumann map plays the role of the inverse
data (see \cite{B08,AM,BM_Sh}). It is worth noting that for
solving the dynamic inverse problem, one can utilize real-valued
boundary controls (see \cite{MM1,MM4}).

\subsection{Forward problem, operators of the Boundary Control method.}
We choose a positive integer $T$ and define the \emph{outer space}
of  the system (\ref{Jacobi_dyn}) as the space of controls or
inputs $\mathcal{F}^T = \mathbb{C}^T$, where $f$ is an element of
$\mathcal{F}^T$. In other words, $f$ is a vector with $T$
components, and each component is a complex number. The inner
product in this space is defined as $(f,g)_{\mathcal{F}^T} =
\sum_{n=0}^{T-1} \overline{f_i} \cdot g_i$, where $\overline{f_i}$
is the complex conjugate of $f_i$. Using this definition, we can
derive a discrete version of the d'Alembert integral
representation formula for the solution of the system
(\ref{Jacobi_dyn}).
\begin{lemma}
The solution to equation (\ref{Jacobi_dyn}) can be represented as follows
\begin{equation}
\label{Jac_sol_rep} u^f_{n,\,t}=\prod_{k=0}^{n-1}
a_kf_{t-n}+\sum_{s=n}^{t-1}w_{n,\,s}f_{t-s-1},\quad n,t\in
\mathbb{N},
\end{equation}
where $w_{n,s}$ satisfies the Goursat problem
\begin{equation*}
\label{Goursat}
\begin{cases}
w_{n,\,s+1}+w_{n,\,s-1}-a_nw_{n+1,\,s}-a_{n-1}w_{n-1,\,s}-b_nw_{n,\,s}=\\
=-\delta_{s,n}(1-a_n^2)\prod_{k=0}^{n-1}a_k,\,n,s\in \mathbb{N}, \,\,s>n,\\
w_{n,\,n}-b_n\prod_{k=0}^{n-1}a_k-a_{n-1}w_{n-1,\,n-1}=0,\quad n\in \mathbb{N},\\
w_{0,\,t}=0,\quad t\in \mathbb{N}_0.
\end{cases}
\end{equation*}
\end{lemma}
Let $\mathcal{F}^\infty=\left\{(f_0,f_1,\ldots)\,|\, f_i\in
\mathbb{C},\, i=0,1,\ldots  \right\}$.
\begin{definition}
For $f,g\in \mathcal{F}^\infty$ we define the convolution
$c=f*g\in \mathcal{F}^\infty$ by the formula
\begin{equation*}
c_t=\sum_{s=0}^{t}f_sg_{t-s},\quad t\in \mathbb{N}\cup \{0\}.
\end{equation*}
\end{definition}

The input $\longmapsto$ output correspondence in the system
(\ref{Jacobi_dyn}) is realized by a \emph{response operator}:
\begin{definition}
For (\ref{Jacobi_dyn}) the \emph{response operator}
$R^T:\mathcal{F}^T\mapsto \mathbb{C}^T$ is defined by the rule
\begin{equation*}
\left(R^Tf\right)_t=u^f_{1,\,t}, \quad t=1,\ldots,T.
\end{equation*}
\end{definition}
This operator plays the role of inverse data. Taking the special
control $f=\delta=(1,0,0,\ldots)$, we define the \emph{response
vector} by
\begin{equation*}
\left(R^T\delta\right)_t=u^\delta_{1,\,t}=: r_{t-1},\quad
t=1,2,\ldots.
\end{equation*}
Then from (\ref{Jac_sol_rep}) we have that
\begin{eqnarray}
\label{R_def}
\left(R^Tf\right)_t=u^f_{1,\,t}=a_0f_{t-1}+\sum_{s=1}^{t-1}
w_{1,\,s}f_{t-1-s}
\quad t=1,\ldots,T.\\
\notag \left(R^Tf\right)=r*f_{\cdot-1}.
\end{eqnarray}
In other words, the response vector
$r=(r_0,r_1,\ldots,r_{T-1})=(a_0,w_{1,1},w_{1,2},\ldots
w_{1,T-1})$ is a convolution kernel of the response operator.

The inverse problem we deal with consists in recovering a Jacobi
matrix (i.e. the sequences $\{a_1,a_2,\ldots\}$, $\{b_1,
b_2,\ldots \}$) and additional parameter $a_0$ from the response
operator.

We introduce the \emph{inner space} (the space of states) of
dynamical system (\ref{Jacobi_dyn}) $\mathcal{H}^T:=\mathbb{C}^T$,
$h\in \mathcal{H}^T$, $h=(h_1,\ldots, h_T)$ with the product
$h,m\in \mathcal{H}^N$,
$(h,m)_{\mathcal{H}^N}:=\sum_{k=1}^N\overline{h_k}{m_k}$, the
space of states. We regard $u^f_{\cdot,\,T}\in \mathcal{H}^T$ (
see (\ref{Jac_sol_rep})).

The input $\longmapsto$ state correspondence of the system
(\ref{Jacobi_dyn}) is realized by a \emph{control operator}
$W^T:\mathcal{F}^T\mapsto \mathcal{H}^T$, defined by the rule
\begin{equation}
\label{cont_oper} \left(W^Tf\right)_n:=u^f_{n,\,T},\quad
n=1,\ldots,T.
\end{equation}
From (\ref{Jac_sol_rep}) we obtain the representation for $W^T$:
\begin{equation*}
\left(W^Tf\right)_n=u^f_{n,\,T}=\prod_{k=0}^{n-1}
a_kf_{T-n}+\sum_{s=n}^{T-1}w_{n,\,s}f_{T-s-1},\quad n=1,\ldots,T.
\end{equation*}
The concept of \emph{boundary controllability} is crucial when
applying the BC method to dynamical systems: \cite{B07,B17}. A
statement equivalent to the boundary controllability of
(\ref{Jacobi_dyn}) can be expressed as follows.
\begin{lemma}
\label{teor_control} The operator $W^T$ is an isomorphism from
$\mathcal{F}^T$ on $\mathcal{H}^T$.
\end{lemma}
The proof easily follows from the triangular representation of
$W^T$:
\begin{equation}
\label{WT}
W^Tf=\begin{pmatrix} u_{1,\,T}\\
u_{2,\,T}\\
\cdot\\
u_{k,\,T}\\
\cdot\\
u_{T,\,T}
\end{pmatrix}=\begin{pmatrix}
a_0& w_{1,1} & w_{1,2} & \ldots & \ldots & w_{1,T-1}\\
0 & a_0a_1 & w_{2,2} & \ldots & \ldots & w_{2, T-1}\\
\cdot & \cdot & \cdot & \cdot & \cdot & \cdot \\
0 & \ldots & \prod_{j=0}^{k-1}
a_j& w_{k,k} & \ldots & w_{k, T-1}\\
\cdot & \cdot & \cdot & \cdot & \cdot & \cdot \\
0 & 0 & 0 & 0 & \ldots & \prod_{k=0}^{T-1} a_{T-1}
\end{pmatrix}
\begin{pmatrix} f_{T-1}\\
f_{T-2}\\
\cdot\\
f_{T-k-1}\\
\cdot\\
f_{0}
\end{pmatrix}.
\end{equation}

The \emph{connecting operator} $C^T: \mathcal{F}^T\mapsto
\mathcal{F}^T$ associated with the system (\ref{Jacobi_dyn}) is
introduced by the bilinear form: for arbitrary $f,g\in
\mathcal{F}^T$ we define
\begin{equation}
\label{C_T_def} \left(C^T
f,g\right)_{\mathcal{F}^T}=\left(u^f_{\cdot,\,T},
u^g_{\cdot,\,T}\right)_{\mathcal{H}^T}=\left(W^Tf,W^Tg\right)_{\mathcal{H}^T}.
\end{equation}
That is $C^T=\left(W^T\right)^*W^T$. The following fact on the
special representation of the connecting operator in terms of
dynamic inverse data is important in the BC method.
\begin{theorem}
The connecting operator $C^T$ is an isomorphism in
$\mathcal{F}^T$, it admits the representation in terms of inverse
data:
\begin{equation}
\label{C_T_repr} C^T=a_0C^T_{ij},\quad
C^T_{ij}=\sum_{k=0}^{T-\max{i,j}}r_{|i-j|+2k},\quad r_0=a_0.
\end{equation}
\begin{equation*}
C^T=
\begin{pmatrix}
r_0+r_2+\ldots+r_{2T-2} & r_1+\ldots+r_{2T-3} & \ldots &
r_{T-2}+r_{T} &
r_{T-1}\\
r_1+r_3+\ldots+r_{2T-3} & r_0+\ldots+r_{2T-4} & \ldots & \ldots
&r_{T-2}\\
\cdot & \cdot & \cdot & \cdot & \cdot \\
r_{T-3}+r_{T-1}+r_{T+1} &\ldots & r_0+r_2+r_4 & r_1+r_3 & r_2\\
r_{T}+r_{T-2}&\ldots &r_1+r_3&r_0+r_2&r_1 \\
r_{T-1}& r_{T-2}& r_2 & r_1 &r_0
\end{pmatrix}
\end{equation*}
\end{theorem}
\begin{proof}
We outline the proof of this statement. First we observe that
$C^T=\left(W^T\right)^*W^T$, so due to Lemma \ref{teor_control},
$C^T$ is an isomorphism in $\mathcal{F}^T$. For fixed $f,g\in
\mathcal{F}^T$ we introduce the \emph{Blagoveshchenskii function}
by the rule
\begin{equation*}
\psi_{n,\,t}:=\left(u^f_{\cdot,\,n},
u^g_{\cdot,\,t}\right)_{\mathcal{H}^T}=\sum_{k=1}^T
\overline{u^f_{k,\,n}}{u^g_{k,\,t}}.
\end{equation*}
We can show that $\psi_{n,\,t}$ satisfies the following difference
equation:
\begin{eqnarray}
\label{Blag_eqn} &\left\{
\begin{array}l
\psi_{n,\,t+1}+\psi_{n,\,t-1}-\psi_{n+1,\,t}-\psi_{n-1,\,t}=h_{n,\,t},\quad n,t\in \mathbb{N}_0,\\
\psi_{0,\,t}=0,\,\, \psi_{n,\,0}=0,
\end{array}
\right. \\
&\text{where}\quad h_{n,\,t}=a_0\left[{g_t}\overline{(Rf)_n}-\overline{f_n}{(Rg)_t}\right].\notag
\end{eqnarray}
We introduce the set
\begin{multline*}
K(n,t):=\left\{(n,t)\cup \{(n-1,t-1),
(n+1,t-1)\}\cup\{(n-2,t-2),(n,t-2),\right.\\
\left.(n+2,t-2)\}\cup\ldots\cup\{(n-t,0),(n-t+2,0),\ldots,(n+t-2,0),(n+t,0)\}\right\}\\
=\bigcup_{\tau=0}^t\bigcup_{k=0}^\tau
\left(n-\tau+2k,t-\tau\right).
\end{multline*}
Then the solution to (\ref{Blag_eqn}) is given by (see \cite{MM1})
\begin{equation*}
\psi_{n,\,t}=\sum_{(k,\tau)\in K(n,t-1)}h(k,\tau)
\end{equation*}
and that $\psi_{T,\,T}=\left(C^Tf,g\right)_{\mathcal{F}^T}$, so we
have
\begin{equation}
\label{C_T_sol}\left(C^Tf,g\right)_{\mathcal{F}^T}=\sum_{(k,\tau)\in
K(T,T-1)}h(k,\tau).
\end{equation}
We note that in the right hand side of (\ref{C_T_sol}) the
argument $k$ runs from $1$ to $2T-1.$ Extending $f\in
\mathcal{F}^T$, $f=(f_0,\ldots,f_{T-1})$ to $f\in
\mathcal{F}^{2T}$ by:
\begin{equation*}
f_T=0,\quad f_{T+k}=-f_{T-k},\,\, k=1,2,\ldots, T-1,
\end{equation*}
we deduce that $\sum_{k,\tau\in K(T,T-1)}
\overline{f_k}{(R^Tg)_\tau}=0$, so (\ref{C_T_sol}) leads to
\begin{multline*}
\left(C^Tf,g\right)=\sum_{k,\tau\in K(T,T-1)}
{g_\tau}\overline{\left(R^{2T}f\right)_k}={g_0}\overline{\left[\left(R^{2T}f\right)_1+\left(R^{2T}f\right)_3+\ldots+\left(R^{2T}f\right)_{2T-1}\right]}\\
+{g_1}\overline{\left[\left(R^{2T}f\right)_2+\left(R^{2T}f\right)_4+\ldots+\left(R^{2T}f\right)_{2T-2}\right]}+\ldots+{g_{T-1}}\overline{\left(R^{2T}f\right)_{T}}.
\end{multline*}
The latter relation yields the representation (\ref{C_T_repr}).
\end{proof}

Let us make some comments on the dependence of the solution and
inverse data on the boundary control and the coefficients of the
matrix $A$:
\begin{remark}
\label{Rem_finite} The formula (\ref{Jac_sol_rep}) implies a
finite wave propagation speed in (\ref{Jacobi_dyn}):  taking
control $f=\delta=(1,0,0\ldots)$ in (\ref{Jacobi_dyn}) we see that
the wave reaches the point $n=N$ exactly at time $t=T$ and on the
front $u^\delta_{N,N}=\prod_{k=0}^{N-1} a_k$. The formula
(\ref{Jac_sol_rep}) implies also the following dependence of
inverse data on coefficients $\{a_n,b_n\}$: for $M\in \mathbb{N}$
the element $u^f_{1,2M-1}$ depends on
$\left\{a_0,a_1,\ldots,a_{M-1}\right\},$
$\left\{b_1,\ldots,b_{M-1}\right\}$.
\end{remark}

\begin{remark}
\label{Rem1} The response operator $R^{2T}$ (or, what is
equivalent, the response vector $(r_0,r_1,\ldots,r_{2T-2})$) is
determined by $\{a_0,\ldots,a_{T-1}\}$, $\{b_1,\ldots,b_{T-1}\}$.
\end{remark}
This observation leads to the natural statement of the dynamic IP
for the system (\ref{Jacobi_dyn}): given the operator $R^{2T}$, to
recover the numbers $\{a_0, \ldots, a_{T-1}\}$ and $\{b_1, \ldots,
b_{T-1}\}$.

We can see from (\ref{R_def}) that $a_0 = r_0$, so we can
immediately determine $a_0$ based on the knowledge of $R^T$. Since
$a_0$ is a multiplier in (\ref{C_T_repr}), we can assume without
loss of generality that the following holds true.
\begin{agreement}
In what follows we assume $a_0=1$.
\end{agreement}

\subsection{Krein equations}

Let $\alpha,\beta\in \mathbb{R}$ and $y(\lambda)$ be a solution to
a Cauchy problem for the following difference equation (we recall
the Agreement $a_0=1$):
\begin{equation*}
\label{y_special} \left\{
\begin{array}l
a_ky_{k+1}+a_{k-1}y_{k-1}+b_ky_k=\lambda y_k,\\
y_0=\alpha,\,\, y_1=\beta.
\end{array}
\right.
\end{equation*}
We set up the \emph{special control problem}: to find a control
$f^T\in \mathcal{F}^T$ that drives the system (\ref{Jacobi_dyn})
to the prescribed state
$\overline{y^T(\lambda)}:=\left(\overline{y_1(\lambda)},\ldots,\overline{y_T(\lambda)}\right)\in
\mathcal{H}^T$ at $t=T$:
\begin{equation}
\label{Control_probl} W^Tf^T=\overline{y^T(\lambda)},\quad
\left(W^Tf^T\right)_k=\overline{y_k(\lambda)},\quad k=1,\ldots,T.
\end{equation}
Note that due to Lemma \ref{teor_control}, this control problem
has a unique solution
$f^T=\left(W^T\right)^{-1}\overline{y^T(\lambda)}$. We define
$\varkappa^T(\lambda)$ by
\begin{equation}
\label{kappa} \left\{
\begin{array}l
\varkappa^T_{t+1}+\varkappa^T_{t-1}=\lambda\varkappa_t,\quad t=0,\ldots,T,\\
\varkappa^T_{T}=0,\,\, \varkappa^T_{T-1}=1.
\end{array}
\right.
\end{equation}
An important fact is that the control $f^T$ is determined by the
operator $R^{2T}$ as follows:
\begin{theorem}
\label{Theor_Krein} The control
$f^T=\left(W^T\right)^{-1}\overline{y^T(\lambda)}$ solving the
special control problem (\ref{Control_probl}), is a unique
solution to the following Krein-type equation in $\mathcal{F}^T$:
\begin{equation}
\label{C_T_Krein}
C^Tf^T=\beta\overline{\varkappa^T(\lambda)}-\alpha
\left(R^T\right)^*\overline{\varkappa^T(\lambda)}.
\end{equation}
\end{theorem}
\begin{proof}
Let $f^T$ be a solution to (\ref{Control_probl}). It can be
checked using (\ref{kappa}) that for any fixed $g\in
\mathcal{F}^T$ we have that
\begin{equation*}
\label{Kr_1}
u^g_{k,\,T}=\sum_{t=0}^{T-1}\left(u^g_{k,\,t+1}+u^g_{k,\,t-1}-\lambda
u^g_{k,\,t}\right)\varkappa^T_t.
\end{equation*}
Using this observation, we can evaluate the quadratic form:
\begin{multline*}
\left(C^Tf^T,g\right)_{\mathcal{F}^T}=\left(u^f(T),u^g(T)\right)_{\mathcal{\mathcal{F}}^T}=\sum_{k=1}^T
\overline{u^f_{k,\,T}}{u^g_{k\,T}}\\
=\sum_{k=1}^T {y_k(\lambda)}{u^g_{k,\,T}}=\sum_{k=1}^T
{y_k(\lambda)}{\sum_{t=0}^{T-1}\left(u^g_{k,\,t+1}+u^g_{k,\,t-1}-\lambda
u^g_{k,\,t}\right)\varkappa^T_t}\\
=\sum_{t=0}^{T-1}{\varkappa^T_t(\lambda)}\left(\sum_{k=1}^T
\left(a_k{u^g_{k+1,\,t}}{y_k}+a_{k-1}{u^g_{k-1,\,t}}{y_k} +
b_k{u^g_{k,\,t}}{y_k}-{\lambda} {u^g_{k,\,t}}{y_k}\right)\right)\\
=\sum_{t=0}^{T-1}{\varkappa^T_t(\lambda)}\left(\sum_{k=1}^T
\left({u^g_{k,\,t}}(a_k{y_{k+1}}+a_{k-1}{y_{k-1}}+b_k{y_k}-
{\lambda y_k} \right)+{u^g_{0,\,t}}{y_1}\right.\\
\left.+a_T{u^g_{T+1,\,t}}{y_T}-{u^g_{1,\,t}}{y_0}-
a_T{u^g_{T,\,t}}{y_{T+1}} \right)
=\sum_{t=0}^{T-1}{\varkappa^T_t(\lambda)}\left(\beta
{g_t}-\alpha{\left(R^Tg\right)_t} \right)\\
=\left(\overline{\varkappa^T(\lambda)}, \left[\beta g-\alpha
\left(R^Tg\right)\right]\right)_{\mathcal{F}^T}=\left(\left[\beta\overline{\varkappa^T(\lambda)}
- \alpha
\left(\left(R^T\right)^*\overline{\varkappa^T(\lambda)}\right)\right],
g\right)_{\mathcal{F}^T}.
\end{multline*}
Which completes the proof due to the arbitrariness of $g$.
\end{proof}

The iterative procedure of recovering $a_k,b_k$ form the solutions
of Krein (\ref{C_T_Krein}) equations $f^\tau\in \mathcal{F}^\tau$
for $\tau=1,\ldots,T$ is given in \cite[Section 3]{MM4}.

\subsection{Factorization method}

In this method one uses the fact that the matrix $C^T$ is a
product of a triangular matrix and its conjugate. We rewrite the
operator $W^T$ as $W^T=W_TJ_T$ where
\begin{equation*}
W^Tf=\begin{pmatrix}
a_0 & w_{1,1} & w_{1,2} & \ldots & w_{1,T-1}\\
0 & a_0a_1 & w_{2,2} &  \ldots & w_{2, T-1}\\
\cdot & \cdot & \cdot & \cdot & \cdot \\
0 & \ldots & \prod_{j=1}^{k-1}a_j& \ldots & w_{k, T-1}\\
\cdot & \cdot & \cdot  & \cdot & \cdot \\
0 & 0 & 0  & \ldots & \prod_{j=1}^{T-1}a_j
\end{pmatrix}
\begin{pmatrix}
0 & 0 & 0 & \ldots & 1\\
0 & 0 & 0 & \ldots  & 0\\
\cdot & \cdot & \cdot & \cdot &  \cdot \\
0 & \ldots & 1& 0  & 0\\
\cdot & \cdot & \cdot & \cdot &  \cdot \\
1 & 0 & 0 & 0 &  0
\end{pmatrix}
\begin{pmatrix} f_{0}\\
f_{2}\\
\cdot\\
f_{T-k-1}\\
\cdot\\
f_{T-1}
\end{pmatrix}.
\end{equation*}
The definition (\ref{C_T_def}) and the invertibility of $W^T$ (cf.
Lemma \ref{teor_control}) imply
\begin{equation*}
C^T=\left(W^T\right)^*W^T,\quad \text{or} \quad
\left(\left(W^T\right)^{-1}\right)^*C^T\left(W^T\right)^{-1}=I.
\end{equation*}
The latter equation can be rewritten as
\begin{equation}
\label{C_T_eqn_ker} \left(\left(W_T\right)^{-1}\right)^* C_T\left(
W_T\right)^{-1}=I,\quad C_T=J_TC^TJ_T.
\end{equation}
Where $\left( W_T\right)^{-1}$ and $C_T$ have  forms
\begin{equation*}
\left(W_T\right)^{-1}=\begin{pmatrix}
q_{1,\,1} & q_{1,\,2} & q_{1,\,3}& \ldots & q_{1,\,T} \\
0 & q_{2,\,2} & q_{2,\,3} &\ldots &..\\
\cdot & \cdot & \cdot & q_{T-1,\,T-1} & q_{T-1,T\,} \\
0 &\ldots &\ldots & 0 & q_{T,\,T}
\end{pmatrix}, \quad
C_T=\left\{m_{ij}\right\}_{i,j=1}^T.
\end{equation*}

We multiply the $k-$th row of $W_T$ by $k-$th column of $\left(
W_T\right)^{-1}$ to get $q_{k,\,k}a_0a_1\ldots a_{k-1}=1$, so
diagonal elements of $\left(W_T\right)^{-1}$ satisfy the relation:
\begin{equation}
\label{AK_form} q_{k,\,k}=\left(\prod_{j=0}^{k-1}a_j\right)^{-1}.
\end{equation}
Multiplying the $k-$th row of $W_T$ by $k+1-$th column of
$\left(W_T\right)^{-1}$ leads to the relation
\begin{equation*}
q_{k,\,k+1}a_0a_1\ldots a_{k-1}+q_{k+1,\,k+1}w_{k,\,k}=0,
\end{equation*}
from where we deduce that
\begin{equation}
\label{AK1_form}
q_{k,\,k+1}=-\left(\prod_{j=0}^{k}a_j\right)^{-2}a_k w_{k,\,k}.
\end{equation}

All aforesaid leads to the equivalent form of (\ref{C_T_eqn_ker}):
\begin{equation}
\label{matr_eq}
\begin{pmatrix}
q_{1,\,1} & 0 & . &  0 \\
q_{1,\,2} & a_{2,\,2} & 0  &.\\
\cdot & \cdot & \cdot & \cdot  \\
q_{1,\,T} & . & .  & q_{T,\,T}
\end{pmatrix}
\begin{pmatrix}
m_{11} & .. & .. &  m_{1T} \\
.. & .. & ..  &..\\
\cdot & \cdot & \cdot & \cdot  \\
m_{T1} &.. &   & m_{TT}
\end{pmatrix}\begin{pmatrix}
q_{1,1} & q_{1,2} & ..& q_{1,T} \\
0 & q_{2,2} & ..  & q_{2,T}\\
\cdot & \cdot & \cdot & \cdot  \\
0 &\ldots &\ldots  & q_{T,T}
\end{pmatrix}=I.
\end{equation}
In the above equality $m_{ij}$ are
entries of $C_T$, the entries $q_{ij}$ are unknown. A direct
consequence of (\ref{matr_eq}) is an equality for determinants:
\begin{equation*}
\det{\left(\left(W_T\right)^{-1}\right)^*}\det{
C_T}\det{\left(W_T\right)^{-1}}=1,
\end{equation*}
which yields
\begin{equation*}
q_{1,\,1}\ldots q_{T,\,T}=\left(\det{ C_T}\right)^{-\frac{1}{2}}.
\end{equation*}
From the above equality we derive that
\begin{equation*}
q_{1,\,1}=\left(\det{C_1}\right)^{-\frac{1}{2}},\quad
q_{2,\,2}=\left(\frac{\det{ C_2}}{\det{
C_1}}\right)^{-\frac{1}{2}},\quad
q_{k,\,k}=\left(\frac{\det{C_k}}{\det{
C_{k-1}}}\right)^{-\frac{1}{2}}.
\end{equation*}
Combining latter relations with (\ref{AK_form}), we deduce that

\begin{equation}
\label{AK} a_k=\frac{\left(\det{
C_{k+1}}\right)^{\frac{1}{2}}\left(\det{
C_{k-1}}\right)^{\frac{1}{2}}}{\det{ C_{k}}}, \quad k=1,\ldots,
T-1.
\end{equation}
Here we set $\det{C_0}=1,$ $\det{C_{-1}=1}$.

Now using (\ref{matr_eq}) we write down the equation on the last
column of $\left(W_T\right)^{-1}$:
\begin{equation*}
\begin{pmatrix}
q_{1,\,1} & 0 & . &  0 \\
q_{1,\,2} & q_{2,\,2} & 0  &.\\
\cdot & \cdot & \cdot & \cdot  \\
q_{1,\,T-1} & . & .  & q_{T-1,\,T-1}
\end{pmatrix}
\begin{pmatrix}
m_{1,\,1} & .. & .. &  m_{1,\,T} \\
.. & .. & ..  &..\\
\cdot & \cdot & \cdot & \cdot  \\
m_{T-1,\,1} &.. &   & m_{T-1,\,T-1}
\end{pmatrix}\begin{pmatrix}
q_{1,\,T} \\
q_{2,\,T}\\
\cdot  \\
q_{T,\,T}
\end{pmatrix}=\begin{pmatrix}
0 \\
0\\
\cdot  \\
0
\end{pmatrix}.
\end{equation*}
As we know $q_{T,\,T}$, we rewrite the above equality in the form
of equation on $(q_{1,\,T},\ldots,q_{T-1,\,T})^*$:
\begin{equation}
\label{matr_eq2}
\begin{pmatrix}
m_{1,\,1} & .. & .. &  m_{1,\,T} \\
.. & .. & ..  &..\\
\cdot & \cdot & \cdot & \cdot  \\
m_{T-1,\,1} &.. &   & m_{T-1,\,T-1}
\end{pmatrix}\begin{pmatrix}
q_{1,\,T} \\
q_{2,\,T}\\
\cdot  \\
q_{T-1,\,T}
\end{pmatrix}=-q_{T,\,T}\begin{pmatrix}
q_{1,\,T} \\
q_{2,\,T}\\
\cdot  \\
q_{T-1,\,T}
\end{pmatrix}.
\end{equation}
Introduce the notation:
\begin{equation*}
C_{k,\,k-1}:=\begin{pmatrix}
m_{1,\,1} & .. & .. &  m_{1,\,k-2} &  m_{1,\,k}\\
.. & .. & ..  &..\\
\cdot & \cdot & \cdot & \cdot  \\
m_{k-1,\,1} &.. &   & m_{k-1,\,k-2}& m_{k-1,\,k}
\end{pmatrix},
\end{equation*}
that is $C_{k,\,k-1}$ is constructed from the matrix $C_{k-1}$ by
substituting the last column by $\begin{pmatrix}m_{1,\,k}\\
\ldots \\ m_{k-1,\,k}\end{pmatrix}$. Then from (\ref{matr_eq2}) we
deduce that:
\begin{equation}
\label{a1}q_{T-1,\,T}=-q_{T,\,T}\frac{\det{
C_{T,\,T-1}}}{\det{C_{T-1}}},
\end{equation}
where we assumed that $\det{C_{0,\,-1}}=0.$ On the other hand,
from (\ref{AK_form}), (\ref{AK1_form}) we see that
\begin{equation}
\label{a2}
q_{T-1,\,T}=\left(\prod_{j=0}^{T-1}a_j\right)^{-1}\sum_{k=1}^{T-1}b_k.
\end{equation}
Equating (\ref{a1}) and (\ref{a2}) gives that
\begin{equation*}
\sum_{k=1}^{T-1}b_k=-\frac{\det{C_{T,\,T-1}}}{\det{
C_{T-1}}},\quad \sum_{k=1}^{T}b_k=-\frac{\det{
C_{T+1,\,T}}}{\det{C_{T}}},
\end{equation*}
from where
\begin{equation}
\label{BK} b_k=-\frac{\det{C_{k+1,k}}}{\det{
C_{k}}}+\frac{\det{C_{k,\,k-1}}}{\det{C_{k-1}}}, \quad k=1,\ldots,
T-1.
\end{equation}

\subsection{Characterization of inverse data.}

Here, we pose the following question: how can we determine if a
given vector $(r_0, r_1, r_2, \ldots, r_{2T-2})$ is a response
vector for a dynamical system described by equation
(\ref{Jacobi_dyn}) with some coefficients $(a_0, \ldots, a_{T-1})$
and $(b_1, \ldots, b_{T-1})$? The answer is given in the following
theorem.
\begin{theorem}
\label{Th_char} The  vector $(r_0,r_1,r_2,\ldots,r_{2T-2})$ is a
response vector for the dynamical system (\ref{Jacobi_dyn}) if and
only if the matrix $C^T$ defined by (\ref{C_T_repr}) is positive
definite.
\end{theorem}
The proof is given in \cite[Theorem 3.2]{MM4}.

\subsection{Discrete Schr\"odinger operator}
The special case of the dynamical system (\ref{Jacobi_dyn}) is
particularly interesting. We assume that the Jacobi matrix
corresponds to a discrete Schr\"odinger operator, meaning that
$a_k = 1$ for all $k \in \mathbb{N} \cup \{0\}$, as described in
\cite{MM1}. In this particular situation, the control operator
(\ref{WT}) has the property that all diagonal elements of the
matrix in (\ref{WT}) are equal to $1$. Thus the determinant of
$W^T$ is also equal to $1$. Because of this fact, the connecting
operator (\ref{C_T_def}) and (\ref{C_T_repr}) has a remarkable
property: the determinant of $C^k$ is always equal to 1 for $k =
1, \ldots, T$. This means that not all elements in the response
vector are independent. Specifically, $r_{2m}$ depends on
$r_{2l+1}$ for $l = 0, \ldots, m-1$. Moreover, this property
characterizes the dynamic data of discrete Schr\"odinger
operators.
\begin{theorem}
The  vector $(1,r_1,r_2,\ldots,r_{2T-2})$ is a
response vector for the dynamical system (\ref{Jacobi_dyn}) with
$a_k=1$ if and only if the matrix $C^T$ (\ref{C_T_repr}) is
positive definite and $\det C^l=1,$ $l=1,\ldots,T$.
\end{theorem}
The proof is provided in \cite[Theorem 5]{MM1}.

We note that for the discrete Schr\"odinger operator there is
another option for recovering unknown coefficients: the equality
(\ref{matr_eq}) is equivalent to the Gelfand-Levitan type
equations, the details are given in \cite[Section 4.3]{MM1}.

\subsection{Inverse problem for complex Jacobi matrices}
\label{Section_Complex_matr}

In \cite{MM18,MM19} the authors solved the inverse dynamic problem
for the complex matrix $A$: that is we assume the entries of
(\ref{Jac_matr}) are complex numbers $\{a_1,\ldots\}$, $\{b_1,
b_2,\ldots  \}$, $a_i\not= 0$, and for this numbers and additional
parameter $0\not=a_0\in \mathbb{C}$ we consider a dynamical system
(\ref{Jacobi_dyn}) with complex boundary control
$f=\left(f_0,f_1,\ldots\right)\in \mathcal{F}^T$. The response
$R^T$ operator for this system is introduced by the same rule
(\ref{resp_oper}).

The IP  as in the real-valued case consists in finding $a_k,$
$b_k$ $k=1,2,\ldots, T-1$ and $a_0$ from the knowledge of response
operator (\ref{resp_oper}) $R^{2T}$.

Note that the operator associated with the matrix $A$ is
nonselfadjoint. In \cite{MM18,MM19} the authors used the
nonselfadjoint variant \cite{AB,BH} of the BC method \cite{B07}.
To do so we introduce the auxiliary dynamical system corresponding
to the adjoint system, which is related to complex conjugate
matrix $A^*=\overline A$:
\begin{equation}
\label{Jacobi_dyn_aux}
\begin{cases}
v_{n,t+1}+v_{n,t-1}-\overline{a_{n}}v_{n+1,t}-\overline{a_{n-1}}v_{n-1,t}-\overline{b_n}v_{n,t}=0,\quad n,t\in \mathbb{N},\\
v_{n,-1}=v_{n,0}=0,\quad n\in \mathbb{N}, \\
v_{0,t}=f_t,\quad t\in \mathbb{N}\cup\{0\}.
\end{cases}
\end{equation}
Operators corresponding to the system (\ref{Jacobi_dyn_aux}) are
equipped with the symbol $\#$. The control operators associated
with (\ref{Jacobi_dyn}) and with \ref{Jacobi_dyn_aux} are
introduced by the same rule (\ref{cont_oper}).
\begin{lemma} The control and response
operators of the system $(\ref{Jacobi_dyn_aux})$ are related with control and response
operators of the original system (\ref{Jacobi_dyn}) by the rule
\begin{equation}
\label{adj_prop} W^T_{\#}=\overline{W^T},\quad
R^T_{\#}=\overline{R^T},
\end{equation}
that is the matrix of $W^T_{\#}$ and response vector $r_{\#}$ are
complex conjugate to the matrix of $W^T$ and vector $r$.
\end{lemma}

For the systems (\ref{Jacobi_dyn}), (\ref{Jacobi_dyn_aux}) we
introduce the \emph{connecting operator} $C^T:
\mathcal{F}^T\mapsto \mathcal{F}^T$ by the sesquilinear form: for
arbitrary $f,g\in \mathcal{F}^T$ we define
\begin{equation*}
\label{C_T_def_compl} \left(C^T
f,g\right)_{\mathcal{F}^T}=\left(u^f_{\cdot,T},
v^g_{\cdot,T}\right)_{\mathcal{H}^T}=\left(W^Tf,W^T_\#g\right)_{\mathcal{H}^T}.
\end{equation*}
\begin{lemma}
The matrix of the connecting operator $C^T$ admits the
representation (\ref{C_T_repr}) in terms of inverse data.
\end{lemma}
The proof is given in \cite[Theorem 1]{MM18}. The relations
(\ref{adj_prop}) imply the following
\begin{remark}
The connecting operator is complex symmetric:
\begin{equation*}
\left(C^T\right)^*=\overline{C^T},
\end{equation*}
\end{remark}

To solve the IP for equation (\ref{Jacobi_dyn}) with a complex
Jacobi matrix, the authors used the same methods as in the case of
real-valued matrices. However, there is a key difference:
\begin{remark}

In the complex-valued case, using Krein equations and the
factorization method enables us to recover the coefficients $a_0$,
$(a_k)^2$, $b_k$, where $k = 1, \ldots$.
\end{remark}
Actually, the impossibility to retrieve $a_k$ is not a flaw of the
method but a characteristic of the problem. We can treat
$u^f_{n,t}$ and response vector $r$ as functions depending on
parameters $a_1,a_2,a_3,\ldots$. Then the following dependence is
valid:
\begin{lemma} Let $u^f_{n,t}$ be a solution to (\ref{Jacobi_dyn}) and let $r$ be the response vector. Then:
\begin{itemize}
\item{1)} The value $u^f_{n,t}$ is odd with respect to
$a_1,a_2,\ldots a_{n-1}$ and even with respect to
$a_n,a_{n+1},\ldots$. \item{2)} The response vector $r$ determined
by $(a_1)^2, (a_2)^2,\ldots$
\end{itemize}
\end{lemma}
The proof is given in \cite[Theorem 3]{MM18}.

Another important difference between complex and real-valued
situations is in the characterization of inverse data. In the case
of a complex matrix, this result takes the following form:
\begin{theorem}
The  vector $(r_0,r_1,r_2,\ldots,r_{2T-2})$
is a response vector for the dynamical system (\ref{Jacobi_dyn})
if and only if the complex symmetric matrices $C^{T-k}$,
$k=0,1,\ldots,T-1$ constructed by (\ref{C_T_repr}) are isomorphism
in $\mathcal{F}^{T-k}$.
\end{theorem}

Below we provide a simple example of the importance of the
condition that the each block $C^{T-k}$, $k=0,1,\ldots T-1$ is an
isomorphism (not only $C^T$ as in Theorem \ref{Th_char}). We take
$r_0=1,$ $r_1=1,$ $r_2=0,$ $r_3=0,$ $r_4=-1$, in this case in
accordance with (\ref{C_T_repr})
\begin{equation*}
C_T=\begin{pmatrix} 0 & 1 &0 \\
1 & 1 & 1\\
0 & 1 & 1
\end{pmatrix},
\end{equation*}
so $C_T$ is an isomorphism, but $C_{T-1}$ is not invertible and
formulas corresponding in complex case to (\ref{AK}) and
(\ref{BK}) are not applicable.

\section{Classical moment problems} \label{Moment_problem}

In this section, we explain how to solve classical moment problems
using the BC method for discrete dynamical systems associated with
the matrices $A$ and $A^N$. First, we describe the problem setup.
We deal with a dynamical system (\ref{Jacobi_dyn}) and an
additional system associated with the finite matrix $A^N$. For all
operators of the BC method related to $A$ and $A^N$ spectral
representation formulas are derived. Next, we focus on the
truncated moment problem, where only a finite number of moments
are known. It is important to note that solving the truncated
moment problem is straightforward and based on solving a
generalized spectral problem. Finally, we outline results on the
existence and uniqueness of various moment problems. In this part
of the section, we mainly use two key elements: the fact that
moments are connected to the entries of the response vector
(\ref{R_def})  and the BC method for dynamical systems with finite
matrices, which allows us to recover spectral data from dynamic
ones. For more information on applying these ideas to different
inverse problems, see \cite{B2001,AMM,AMMP}.

\subsection{Setup of the problem}

\noindent{\bf \newline Hamburger moment problem. } Given a
sequence of real numbers $s_0,s_1,s_2,\ldots$ find a Borel measure
$d\rho(\lambda)$ on $\mathbb{R}$ such that equalities
(\ref{Moment_eq}) hold.

\noindent{\bf Stieltjes moment problem. } Given a sequence of real
numbers $s_0,s_1,s_2,\ldots$ find a Borel measure $d\rho(\lambda)$
with support on the half-line $\mathbb{R}_+=[0,+\infty )$ such
that $d\rho(\lambda)$ is a solution to the Hamburger moment
problem.

\noindent{\bf Hausdorff moment problem. } Given a sequence of real
numbers $s_0,s_1,s_2,\ldots$ find a Borel measure $d\rho(\lambda)$
with support on the interval $[0,1]$ such that $d\rho(\lambda)$ is
a solution to the Hamburger moment problem.

The moment problem (moment sequence $s_0, s_1,\ldots$) is called
\emph{determinate} if a solution exists and is unique, if a
solution is not unique, it is called \emph{indeterminate}, the
same notations are used in respect to corresponding measures.

The set of solutions of Hamburger moment problem is denoted by
$\mathcal{M}_H$.

With a set of moments we associate the following Hankel matrices:
\begin{equation}
\label{Hankel}
S^N_{m}:=\begin{pmatrix} s_{2N-2+m} & s_{2N-3+m} & \ldots & s_{N-1+m}\\
s_{2N-3+m} & \ldots & \ldots & \ldots\\
\cdot & \cdot & \ldots & s_{1+m} \\
s_{N-1+m} & \ldots & s_{1+m} & s_{m}
\end{pmatrix},\quad m=0,1,\dots,
\end{equation}
When $m = 0$, we write $S^N$ as $S^N_0$. Note that classical
Hankel matrices usually associated with moments are
$S_N=J_NS^NJ_N$. The main questions related to moment problems
concern the existence and uniqueness of a measure that satisfies
the equation (\ref{Moment_eq}) \cite{AH,Schm,S}. To solve these
problems, we utilize the dynamic IP for both the system
(\ref{Jacobi_dyn}) and similar finite-dimensional systems.

\subsection{Spectral representation of connecting and response
operators}

\subsubsection{Finite Jacobi matrix} \label{Oper_fin}

Consider a dynamical system with discrete time associated with a
finite real Jacobi matrix $A^N$:
\begin{equation}
\label{Jacobi_dyn_int} \left\{
\begin{array}l
v_{n,t+1}+v_{n,t-1}-a_nv_{n+1,t}-a_{n-1}v_{n-1,t}-b_nv_{n,t}=0,\,\, t\in \mathbb{N}_0,\,\, n\in 1,\ldots, N,\\
v_{n,\,-1}=v_{n,\,0}=0,\quad n=1,2,\ldots,N+1, \\
v_{0,\,t}=f_t,\quad v_{N+1,\,t}=0,\quad t\in \mathbb{N}_0,
\end{array}\right.
\end{equation}
we treat the real sequence $f=(f_0,f_1,\ldots)$ as a
\emph{boundary control}. Fixing a positive integer $T$, we take
control $f\in \mathcal{F}^T$.

The solution to the system (\ref{Jacobi_dyn_int}) is denoted by
$v^f$. Note that (\ref{Jacobi_dyn_int}) is a discrete analog of an
initial boundary value problem on an interval for the wave
equation with a potential. At the left end of the interval, the
Dirichlet control is applied, and at the right end, the Dirichlet
boundary condition is imposed.

Introduce the operator $A^N:\mathcal{H}^N\mapsto \mathcal{H}^N$
\begin{equation*}
(A^N\psi)_n=\begin{cases}b_1\psi_1+a_1\psi_2,\quad n=1,\\
a_{n}\psi_{n+1}+a_{n-1}\psi_{n-1}+b_n\psi_n,\quad
2\leqslant n\leqslant N-1,\\
a_{N-1}\psi_{N-1}+b_N\psi_N,\quad n=N.
\end{cases}
\end{equation*}
Denote by $\phi=\left\{\phi_n\right\},$ $n=0,1,2,\ldots$ the
solution to the Cauchy problem for the following difference
equation
\begin{equation}
\label{Phi_def}
\begin{cases} a_n\phi_{n+1}+a_{n-1}\phi_{n-1}+b_n\phi_n=\lambda\phi_n,\quad n\geqslant 1,\\
\phi_0=0,\,\,\phi_1=1,
\end{cases}
\end{equation}
where $\lambda\in \mathbb{C}$. Each $\phi_n$ is a polynomial of
$\lambda$ of degree $n-1$. Denote by $\{\lambda_k\}_{k=1}^N$ the
roots of the equation $\phi_{N+1}(\lambda)=0$, it is known
\cite{AH,S} that they are real and distinct. We introduce the
vectors $\phi(\lambda),\phi^k\in \mathbb{R}^N$ by the rules
\begin{equation}
\label{Phi_k}
\phi(\lambda):=\begin{pmatrix}\phi_1(\lambda)\\
\phi_2(\lambda)\\ \cdot\\ \phi_N(\lambda)\end{pmatrix},\quad \phi^k:=\begin{pmatrix}\phi_1(\lambda_k)\\
\phi_2(\lambda_k)\\ \cdot\\ \phi_N(\lambda_k)\end{pmatrix},\quad
k=1,\ldots,N,
\end{equation}
and define the numbers $\omega_k$ by
\begin{equation}
\label{Ortog} (\phi^k,\phi^l)=\delta_{kl}\omega_k,\quad
k,l=1,\ldots,N,
\end{equation}
where $(\cdot,\cdot)$ is a scalar product in $\mathbb{C}^N$.
Denote by $\mathcal{T}_k(2\lambda)$ the Chebyshev polynomials of
the second kind: they are obtained as a solution to the following
Cauchy problem:
\begin{equation*}
\label{Chebysh} \left\{
\begin{array}l
\mathcal{T}_{t+1}+\mathcal{T}_{t-1}-\lambda \mathcal{T}_{t}=0,\\
\mathcal{T}_{0}=0,\,\, \mathcal{T}_1=1.
\end{array}
\right.
\end{equation*}
In \cite{MM1,MM4} the following representation formula for the
solution of (\ref{Jacobi_dyn_int}) was proved:
\begin{proposition}
The solution to (\ref{Jacobi_dyn_int}) admits the representation
\begin{equation}
\label{Jac_sol_rep_int} v^f_{n,t}=
\begin{cases}
\sum_{k=1}^N c_t^k\phi^k_n,\quad n=1,\ldots,N,\\
f_t,\quad n=0.
\end{cases},\quad
c^k=\frac{1}{\omega_k}\mathcal{T}\left(\lambda_k\right)*f.
\end{equation}
\end{proposition}
The proof is given in \cite[Proposition 1]{MM4}.

The \emph{inner space} of dynamical system (\ref{Jacobi_dyn_int})
is $\mathcal{H}^N:=\mathbb{C}^N$, its elements are complex
columns. By (\ref{Jac_sol_rep_int}) we have that
$v^f_{\cdot,\,T}\in \mathcal{H}^N$. For the system
(\ref{Jacobi_dyn_int}) the \emph{control operator}
$W^T_{N}:\mathcal{F}^T\mapsto \mathcal{H}^N$ is defined by the
rule
\begin{equation*}
W^T_{N}f:=v^f_{n,\,T},\quad n=1,\ldots,N.
\end{equation*}

As well as in (\ref{Jacobi_dyn}), the input $\longmapsto$ output
correspondence in the system (\ref{Jacobi_dyn_int}) is realized by
a \emph{response operator}: $R^T_N:\mathcal{F}^T\mapsto
\mathbb{C}^T$, defined by the formula
\begin{equation} \label{R_def_int}
\left(R^T_Nf\right)_t=v^f_{1,\,t}, \quad t=1,\ldots,T.
\end{equation}
This operator has a form of a convolution:
\begin{equation*}
\left(R^T_Nf\right)_t= \sum_{s=0}^{t}r_s^Nf_{t-s-1}\quad
\text{or}\quad R^T_Nf=r^N*f_{\cdot-1},
\end{equation*}
where the convolution kernel is called a \emph{response vector}:
$r^N=(r^N_0,r^N_1,\ldots,r^N_{T-1})$.

The \emph{connecting operator} $C^T_{N}: \mathcal{F}^T\mapsto
\mathcal{F}^T$ for the system (\ref{Jacobi_dyn_int}) is defined
via the sesquilinear form: for arbitrary $f,g\in \mathcal{F}^T$
one has that
\begin{equation*}
\left(C^T_{N} f,g\right)_{\mathcal{F}^T}=\left(v^f_{\cdot,\,T},
v^g_{\cdot,\,T}\right)_{\mathcal{H}^N}=\left(W^T_{N}f,W^T_{N}g\right)_{\mathcal{H}^N},\quad
C^T_N=\left(W^T_{N}\right)^*W^T_{N}.
\end{equation*}

The finiteness of s speed of a wave propagation in the system
(\ref{Jacobi_dyn_int}) leads to the following (cf. Remark
\ref{Rem1}):
\begin{remark}
\label{Rem11} The components of the response vector
$(r_0^N,r_1^N,\ldots,r_{2N-2}^N)$ are determined by the values of
$\{a_0,\ldots,a_{N-1}\}$, $\{b_1,\ldots,b_{N}\}$,  and they do not
depend on the Dirichlet condition at $n=N+1$. However, the
components $r_{2N-1}^N$, $r_{2N}^N$, and so on, do take into
account the boundary condition at $n=N+1$
\end{remark}

The kernel of the response operator (\ref{R_def_int}) is given by
\begin{equation}
\label{con1}
r_{t-1}^{N}=\left(R^T_{N}\delta\right)_t=v^\delta_{1,\,t},\quad
t=1,\ldots.
\end{equation}
The spectral function of operator $A^N$ is introduced by the rule
\begin{equation}
\label{measure}
\rho^{N}(\lambda)=\sum_{\{k\,|\,\lambda_k<\lambda\}}\frac{1}{\omega_k},
\end{equation}
then from (\ref{Jac_sol_rep_int}), (\ref{con1}) we immediately
deduce.
\begin{proposition}
The solution to (\ref{Jacobi_dyn_int}), the response vector of
(\ref{Jacobi_dyn_int}) and entries of the matrix of the connecting
operator $C^T_N$ admit the following spectral representations:
\begin{gather}
v^f_{n,t}=\int_{-\infty}^\infty \sum_{k=1}^t
\mathcal{T}_k(\lambda)f_{t-k}\phi_n(\lambda)\,d\rho^{N}(\lambda),\quad n,t\in \mathbb{N},\notag 
\\
r_{t-1}^{N}=\int_{-\infty}^\infty
\mathcal{T}_t(\lambda)\,d\rho^{N}(\lambda),\quad t\in \mathbb{N},\label{Resp_int_spectr}\\
\{C^T_{N}\}_{l+1,\,m+1}=\int_{-\infty}^\infty
\mathcal{T}_{T-l}(\lambda)\mathcal{T}_{T-m}(\lambda)\,d\rho^{N}(\lambda),
\,\, l,m=0,\ldots,T-1\label{SP_mes_int} .
\end{gather}
\end{proposition}
The proof is given in \cite[Section 4]{MM4}.

\subsubsection{Semi-infinite Jacobi matrix} \label{Section_spectr_repr}

With the matrix $A$ we associate the symmetric operator $A$ in the
space $l_2$ (we keep the same notation), defined on the finite
sequences:
\begin{equation}
\label{domA}
D(A)=\left\{\varkappa=(\varkappa_0,\varkappa_1,\ldots)\,|\,\exists
N(\varkappa)\in \mathbb{N}, \,\varkappa_n=0,\,\text{for
$n\geqslant N(\varkappa)$}\right\},
\end{equation}
and given by the rule
\begin{equation}
\label{A_oper}
\begin{cases}
(A\theta)_1&=b_1\theta_1+a_1\theta_2,\quad n=1,\\
(A\theta)_n&=a_{n}\theta_{n+1}+a_{n-1}\theta_{n-1}+b_n\theta_n,\quad
n\geqslant 2.
\end{cases}
\end{equation}

For a given sequence
$\varkappa=\left(\varkappa_1,\varkappa_2,\ldots\right)$ we define
a new sequence
\begin{equation*}
\left(G\varkappa\right)_n=\begin{cases}
b_1\theta_1+a_1\theta_2,\quad n=1,\\
a_{n}\theta_{n+1}+a_{n-1}\theta_{n-1}+b_n\theta_n,\quad n\geqslant
2.
\end{cases}
\end{equation*}

The adjoint operator $A^*\varkappa:=G\varkappa$ is defined on the
domain
\begin{equation*}
D\left(A^*\right)=\left\{\varkappa=(\varkappa_0,\varkappa_1,\ldots)\in
l_2\,|\, (G\varkappa)\in l_2\right\}.
\end{equation*}

We note that the first line in (\ref{A_oper}) implies a zero
boundary condition at $n=0$, and thus $A$ can have deficiency
indices $(0,0)$ or $(1,1)$. In the limit point at infinity case,
which is when $A$ has the deficiency indices $(0,0)$, one can say
that $A$ is essentially self-adjoint. In the limit circle case,
which is when $A$ has deficiency indices $(1,1)$), we denote by
$p(\lambda)=(p_1(\lambda),p_2(\lambda),\ldots)$,
$q(\lambda)=(q_1(\lambda),q_2(\lambda),\ldots)$ two solutions of
difference equation in  (\ref{Phi_def}) satisfying Cauchy data
$p_1(\lambda)=1,$ $p_2(\lambda)=\frac{\lambda-b_1}{a_1}$,
$q_1(\lambda)=0,$ $q_2(\lambda)=\frac{1}{a_1}$. Then \cite[Lemma
6.22]{Schm}
\begin{equation*}
D\left(A^*\right)=D(A)+\mathbb{C}p(0)+\mathbb{C}q(0).
\end{equation*}

All self-adjoint extensions of $A$ are parameterized by $\alpha\in
\mathbb{R}\cup\{\infty\}$, are denoted by $A_\alpha$ and defined
on the domain
\begin{equation*}
D\left(A_\alpha\right)=\begin{cases} D(A)+\mathbb{C}(q(0)+\alpha p(0)),\quad \alpha\in \mathbb{R}\\
D(A)+\mathbb{C}p(0),\quad \alpha=\infty.
\end{cases}
\end{equation*}

All the details the reader can find in \cite{S,Schm}. By
$[\cdot,\cdot]$ we denote the scalar product in $l_2$. We
introduce the measure
$d\rho_{\alpha}(\lambda)=\left[dE^{A_{\alpha}}_\lambda
e_1,e_1\right]$, where $dE^{A_\alpha}_\lambda$ is the
projection-valued spectral measure of $A_\alpha$ such that
$E^{A_\alpha}_{\lambda-0}=E^{A_\alpha}_{\lambda}$. The results of
\cite{AT} and \cite[Section 5]{S} imply that in the limit circle
case $d\rho^N\to d\rho_{\alpha^*}$ $*-$weakly as $N\to\infty,$
where
\begin{equation}
\label{alpha_def}
\alpha^*=-\lim_{n\to\infty}\frac{q_n(0)}{p_n(0)}.
\end{equation}
In the limit point case $d\rho^N$ converges  $*-$weakly as
$N\to\infty,$ to the spectral measure of $A$ which we also denote
as $d\rho_{\alpha^*}$ for uniformity.

The Remark \ref{Rem11} in particular implies that
\begin{eqnarray}
\label{R_eqv} R^{2N-2} = R^{2N-2}_{N},\\
\label{Rav_bc} u^f_{n,\,t}=v^f_{n,\,t},\quad n\leqslant t\leqslant
N,\quad \text{and}\quad W^N=W^N_{N}.
\end{eqnarray}
Thus due to (\ref{R_eqv}), we have that $r_{t-1}=r_{t-1}^{N}$,
$t=1,\ldots,2N$. On the other hand, taking into account
(\ref{Rav_bc}), we can see that $C^T=C^T_{N}$ with $T\leqslant N$.
Thus, in (\ref{Resp_int_spectr}), (\ref{SP_mes_int}) tending
$N\to\infty$, we obtain
\begin{proposition}
The solution to (\ref{Jacobi_dyn}), the entries of the response
vector of (\ref{Jacobi_dyn}) and of the matrix of the connecting
operator $C^T$ admit the spectral representations:
\begin{gather}
u^f_{n,t}=\int_{-\infty}^\infty \sum_{k=1}^t
\mathcal{T}_k(\lambda)f_{t-k}\phi_n(\lambda)\,d\rho_{\alpha^*}(\lambda),\quad
n,t\in \mathbb{N},\label{Jac_sol_spectr_inf}\\
r_{t-1}=\int_{-\infty}^\infty
\mathcal{T}_t(\lambda)\,d\rho_{\alpha^*}(\lambda),\quad t\in
\mathbb{N},\label{Resp_spectr}\\
\{C^T\}_{l+1,\,m+1}=\int_{-\infty}^\infty
\mathcal{T}_{T-l}(\lambda)\mathcal{T}_{T-m}(\lambda)\,d\rho_{\alpha^*}(\lambda),
\quad l,m=0,\ldots,T-1. \label{SP_mes_d}
\end{gather}
\end{proposition}
The proof can be found in \cite[Section4]{MM4}. Note that in the
limit circle case in (\ref{Jac_sol_spectr_inf}),
(\ref{Resp_spectr}) and (\ref{SP_mes_d}), one can replace
$d\rho_{\alpha^*}$ with $d\rho_\alpha$ for any $\alpha\in
\mathbb{R}\cup\{\infty\}$. Also, (\ref{Resp_spectr}) shows the
relationship between moments and components of the response
vector.
\begin{proposition}
The elements of the response vector are
related to the moments according to the following rule:
\begin{equation}
\label{Resp_and_moments}
\begin{pmatrix}
r_0\\
r_1\\
\ldots \\
r_{n-1}
\end{pmatrix}=\Lambda_n\begin{pmatrix}
s_0\\
s_1\\
\ldots \\
s_{n-1}
\end{pmatrix}.
\end{equation} where the entries of the matrix $\Lambda_n\in \mathbb{M}^{n}$ are
given by
\begin{equation}
\label{LambdaMatr}
\{\Lambda_n\}_{ij}=a_{ij}=\begin{cases} 0,\quad \text{if $i<j$},\\
0,\quad \text{if $i+j$ is odd,}\\
E_{\frac{i+j}{2}}^j(-1)^{\frac{i+j}{2}+j},\quad \text{otherwise}
\end{cases}
\end{equation}
where $E_n^k$ are binomial coefficients.
\end{proposition}
The proof is given in \cite[Proposition 6]{MM6}

\subsection{Truncated moment problem}
\label{Section_tranc}

First we consider the following problem
\begin{definition}
By a solution of a truncated moment problem of order $N$ we call a
Borel measure $d\rho(\lambda)$ on $\mathbb{R}$ such that
equalities (\ref{Moment_eq}) with this measure hold for
$k=0,1,\ldots,2N-2.$
\end{definition}

The spectral representation of the response vector (see equation
(\ref{Resp_spectr})) and the results from the previous sections
lead to the following:
\begin{remark}
The knowledge of a finite set of moments
$\{s_0,s_1,\ldots,s_{2N-2}\}$ is equivalent to knowledge of
$\{r_0,r_1,\ldots,r_{2N-2}\}$, from which we can recover the
Jacobi matrix $A^N \in \mathbb{M}^N$. With this matrix one
associates the dynamical system with Dirichlet boundary condition
at $n=N+1$ given by equation (\ref{Jacobi_dyn_int}), or $N\times
N$ block in the semi-infinite Jacobi matrix in
(\ref{Jacobi_dyn_int}).
\end{remark}

First, we describe the naive method of solving the truncated moment
problem:
\begin{procedure}
\label{Proc1}
\begin{itemize}\item[]
\item[1)]Calculate by (\ref{Resp_and_moments}) the values of $r_0,
r_1, r_2, \ldots, r_{2N-2}$ based  on the given numbers $s_0, s_1,
\ldots, s_{2N-2}$.

\item[2)] Recover the $N \times N$ Jacobi matrix $A^N$ using the
formulas (\ref{AK}) and (\ref{BK}).

\item[3)] Recover the spectral measure for the finite Jacobi
matrix $A^N$, which has a specified Dirichlet or other
self-adjoint condition at  $n=N+1$.

\item[3')] Extend the Jacobi matrix $A^N$ to a finite Jacobi
matrix $A^M$ in an arbitrary way, where $M$ is greater than $N$.
Then, specify an arbitrary self-adjoint condition at the position
$n=M+1$. Finally, reconstruct the spectral measure of $A^M$.

\item[3'')] Extend the Jacobi matrix $A^N$ arbitrarily to form a
semi-infinite Jacobi matrix $A$. Then,  recover the spectral
measure of $A$.
\end{itemize}
\end{procedure}
Every measure obtained in $3)$, $3')$, and $3'')$ provides a
solution to the truncated moment problem.

\begin{remark}
\label{Rem2} In view of the described dependence of $a_k,$ $b_k$
on the components of response vector (see \ref{AK}), (\ref{BK}),
see Remarks \ref{Rem1} and \ref{Rem11}, in the formulas
(\ref{Jac_sol_spectr_inf})--(\ref{SP_mes_d})  one can change
$d\rho_{\alpha^*}(\lambda)$ to the appropriate (depending on the
range of variation of parameters $n,\,t$) measure taken from
Procedure \ref{Proc1}.
\end{remark}

\begin{remark}
\label{Rem3} The fact that dynamic data are related to spectral
one by simple formulas (\ref{Resp_spectr}), (\ref{SP_mes_d}) and
taking in account Remark \ref{Rem2}, allows one to reduce spectral
IP for (\ref{Jacobi_dyn}) or (\ref{Jacobi_dyn_int}) on the
reconstruction $A$ ($A^N$) by the spectral measure to a dynamic IP
on the reconstruction $A$ ($A^N$) by the response operator.
\end{remark}

Note that in the proposed approach, the step of finding the matrix
$A^N$ is resource-intensive and, apparently, redundant: we are
interested in the measure, not the operator. In their works
\cite{MM5,MM6}, the authors proposed an alternative method, which
directly reconstructs the spectral measure associated with the
Jacobian matrix $A^N$, based on moments (actually, the operator
$C^N$) without the need to reconstruct the Jacobian matrix itself.
Here we outline this method.

By $f^k$, where  $k=1,\ldots,N$ we denote the control that drives
the system described by equation (\ref{Jacobi_dyn_int}) towards  a
specific state $\phi^k$,  defined in equation (\ref{Phi_def}).
\begin{equation*}
\label{Control_var} W^Nf_k=\phi^k,\quad k=1,\ldots,N.
\end{equation*}
Due to Lemma \ref{teor_control} and (\ref{Rav_bc}), there exists
a unique control for every $k$.

The controllability of the dynamical system described by equation
(\ref{Jacobi_dyn_int}) (as stated in Lemma \ref{teor_control})
allows us to use the method of variations for minimizing the
quadratic form $U[f]:=(C^Nf,f)_{\mathcal{F}^N}$ in the space
$\mathcal{F}^N$. This method enables us to recover the spectral
data of $A^N$, including the spectrum
$\left\{\lambda_k\right\}_{k=1}^N$, and the coefficients defined
by equation (\ref{Ortog}), from the dynamic data (response
vector), which is encoded in the matrix $C^N$ according to
equation (\ref{C_T_repr}), see also \cite{B2001}.

We introduce the shift operator
\begin{gather*}
V^N: \mathcal{F}^N\mapsto \mathcal{F}^{N},\\
\left(V^Nf\right)_n=f_{n-1},\quad n=1,\ldots,N-1,\quad
\left(V^Nf\right)_0=0,
\end{gather*}
and denote by $E:\mathcal{F}^{N}\mapsto \mathcal{F}^{N+1}$ the
embedding operator: $\left(Ef\right)_{k}=f_k$, $k=0,1,N-1$, and
$\left(Ef\right)_{N}=0$; then the adjoint operator $E^*:
\mathcal{F}^{N+1}\mapsto \mathcal{F}^{N}$ is a projection. We also
introduce
\begin{equation*}
B^N:=E^*\left(V^{N+1}\right)^*C^{N+1}E+C^NV^N.
\end{equation*}
Then corresponding Euler-Lagrnage equations has the following
form:
\begin{theorem}
The spectrum of $A^{N}$ and controls $f_k$,
$k=1,\ldots,N$ are the spectrum and the eigenvectors of the
following generalized spectral problem:
\begin{equation}
\label{m_eqn} B^Nf_k=\lambda_kC^Nf_k,\quad k=1,\ldots,N.
\end{equation}
\end{theorem}

Since (\ref{m_eqn}) is linear, solving this equation one obtains
spectrum $\left\{\lambda_k\right\}_{k=1}^N$ and controls
$\left\{\widetilde f_k\right\}_{k=1}^N$, such that $W^N\widetilde
f_k=\beta_k\phi^k$ for some $\beta_k\in
\mathbb{R}\backslash\{0\}$, $k=1,\ldots,N$. Then the measure of
operator $A^N$ with Dirichlet boundary condition at $n=N+1$ can be
recovered by the following

\begin{procedure}
\begin{itemize}
\item[]
\item[1)] Normalize the controls $\widetilde f_k$ by setting the condition
$\left(C^N\widetilde f_k,\widetilde f_k\right)_{\mathcal{F}^N}=1$, where
$k=1,\ldots,N$.

\item[2)] Observe that $W^N\widetilde f^k=\alpha_k\phi^k$ for some
$\alpha_k\in \mathbb{R}\backslash\{0\}$, where the constant is
defined by $\alpha_k=\alpha_k\phi^k_1=\left(W^N\widetilde
f_k\right)_1=\left(R\widetilde f_k\right)_N$, $k=1,\ldots,N$.

\item[3)] Coefficients (\ref{Ortog}) are given by
$\omega_k={\alpha_k^2}$, $k=1,\ldots,N$.

\item[4)] Recover the measure by (\ref{measure}).

\end{itemize}
\end{procedure}
In the moment problems theory, it is more natural to formulate
results  in terms of Hankel matrices (\ref{Hankel}). To do so we
define the matrix
\begin{equation*}
\widetilde\Lambda_N:=J_N\Lambda_N J_N.
\end{equation*}
\begin{proposition}
\label{CT_ST_prop} The following relation holds:
\begin{eqnarray*}
C^N=\widetilde\Lambda_N S_0^N
\left(\widetilde\Lambda_N\right)^*,\\
B^N =\widetilde\Lambda_N S_1^N
\left(\widetilde\Lambda_N\right)^*.
\end{eqnarray*}
The generalized spectral problem (\ref{m_eqn}) upon introducing
the notation $g_k=\left(\widetilde\Lambda_N\right)^*f_k$ is
equivalent to the following generalized spectral problem:
\begin{equation*}
S^N_1g_k=\lambda_k S^N_0g_k.
\end{equation*}
\end{proposition}

\subsection{Hamburger, Stieltjes and Hausdorff moment problems, existence and uniqueness. }

We consider two special solutions for the equation
(\ref{Phi_def}). The  first one is denoted as $\varphi(\lambda)$
and is obtained by setting $\phi_0=0$ and $\phi_1=1$. The second
solution is denoted as $\xi(\lambda)$ and corresponds to the
Cauchy data $\xi_0=-1$ and $\xi_1=0$.

It is well-known that the uniqueness of the solution to a moment
problem is related to the deficiency indices of the operator $A$
(see \cite{AH,S} ). We provide here the well-known results on the
discrete version of Weyl's limit point-circle theory. This theory
provides an answer to the question about the index of $A$.
\begin{proposition}
\label{Hamburher_crt}The limit circle case holds, meaning that the
operator $A$ has indices $(1,1)$ if and only if one of the
following conditions are fulfilled:
\begin{itemize}
\item $\varphi(\lambda),\xi(\lambda)\in l^2$ for some
$\lambda\in \mathbb{R}$,

\item $\varphi(\lambda),\varphi'(\lambda)\in l^2$ for some
$\lambda\in \mathbb{R}$,

\item $\xi(\lambda),\xi'(\lambda)\in l^2$ for some $\lambda\in
\mathbb{R}$.
\end{itemize}
\end{proposition}

Using Proposition \ref{Hamburher_crt} and the Krein equations
(\ref{C_T_Krein}), as well as the relationship between the
elements of the response vector and the moments
(\ref{Resp_spectr}), the authors in \cite{MM5, MM6} established
the following results:

\begin{theorem}
\label{Th_char_new} The  set of numbers $(s_0,s_1,s_2,\ldots)$ are
the moments of a spectral measure corresponding to the Jacobi
operator $A$ if and only if
\begin{equation*}
\text{the matrix $S^N_0$ is positive definite for
all $N\in \mathbb{N}$.}
\end{equation*}
The Hamburger moment problem is indeterminate if and only if
\begin{equation*}
\lim_{N\to\infty}\left(\left(C^N\right)^{-1}\Gamma_N,\Gamma_N\right)_{\mathcal{F}^N}<+\infty,\quad
\lim_{N\to\infty}\left(\left(C^N\right)^{-1}\Delta_N,\Delta_N\right)_{\mathcal{F}^N}<+\infty,
\end{equation*}
where
\begin{equation*}
\Gamma_N:=\begin{pmatrix} \mathcal{T}_N(0)\\
\mathcal{T}_{N-1}(0)\\\ldots\\
\mathcal{T}_1(0)\end{pmatrix},\quad \Delta_N=\begin{pmatrix} \mathcal{T}_N'(0)\\
\mathcal{T}_{N-1}'(0)\\\ldots\\
\mathcal{T}_1'(0)\end{pmatrix}.
\end{equation*}
\end{theorem}

\begin{theorem}
\label{Th_char_new1} The  set of numbers $(s_0,s_1,s_2,\ldots)$
are the moments of a spectral measure, supported on $(0,+\infty)$,
corresponding to Jacobi operator $A$ if and only if
\begin{equation*}
\text{matrices $S^N_0$ and $S^N_1$ are positive definite for all
$N\in \mathbb{N}$.}
\end{equation*}
The Stieltjes moment problem is indeterminate if and only if the
following relations hold:
\begin{eqnarray*}
M_\infty=\lim_{N\to\infty}\left(\left(C^N\right)^{-1}\Gamma_N,\Gamma_N\right)_{\mathcal{F}^N}<+\infty,\\
L_\infty=\lim_{N\to\infty}\frac{\left(\left(C^{N}\right)^{-1}\left(R^{N}\right)^*\Gamma_{N},e_1\right)}
{\left(\left(C^{N}\right)^{-1}\Gamma_{N},e_1\right)}<+\infty.
\end{eqnarray*}
\end{theorem}
In Section \ref{Section_KreinString} we consider the dynamic IP
for a Krein string, example of which is (\ref{DnSst}) with the
infinite Jacobi matrix of a special form. Note that the conditions
on $M_{\infty}$ and $L_{\infty}$ mentioned above are understood as
that the mass and length of a string are finite, as explained in
references \cite{AH,MM6}.

\begin{theorem}
\label{Th_char_new2} The  set of numbers $(s_0,s_1,s_2,\ldots)$
are the moments of a spectral measure, supported on $[0,1]$,
corresponding to operator $A$ if and only if the condition
\begin{equation*}
S^N_0\geqslant S^N_1 > 0\quad \text{holds for all $N\in
\mathbb{N}$}
\end{equation*}
The Hausdorff moment problem is determinate.
\end{theorem}

The proofs of Theorems \ref{Th_char_new}, \ref{Th_char_new1},
\ref{Th_char_new2} are given in \cite[Proposition 7,8,9]{MM6}.

\section{Moment problem and the Toda lattices}

In this section, we outline the method of constructing the
solution  to (\ref{Toda_eq}) and (\ref{Toda_init}) for a wide
class of unbounded sequences (\ref{Toda_init}).

With the system (\ref{Toda_eq}) we associate two operators acting
in  $l^\infty$, with domains
\begin{equation*}
D(A(t))=D(P(t))=\left\{\varkappa=(\varkappa_0,\varkappa_1,\ldots)\,|\,\exists
N(\varkappa)\in \mathbb{N}, \,\varkappa_n=0,\,\text{for
$n\geqslant N(\varkappa)$}\right\}.
\end{equation*}
The first one, which we call $A(t)$, is given by a self-adjoint
extension (see Section \ref{Section_spectr_repr}) of a symmetric
operator represented by a matrix given by (\ref{Jac_matr}). The
entries of this matrix depend on the variable $t$. $P(t)$ is an
extension of a skew-symmetric operator defined by the following
rule.
\begin{align*}
&\left(P(t)f\right)_1=a_1(t)f_2,\\
&\left(P(t)f\right)_n=a_n(t)f_{n+1}-a_{n-1}(t)f_{n-1},\quad
n=2,\ldots.
\end{align*}

It is a well-known fact \cite{To,Te} that the system
(\ref{Toda_eq}) is equivalent to the following operator equation:
\begin{equation}
\label{Toda_eqviv} \frac{dA}{dt}=PA-AP,
\end{equation}
at least for the bounded coefficients $a_k(t),b_k(t)$
$k=1,2,\ldots$.

Let $d\rho^t(\lambda)$ denote the spectral measure of $A(t)$.
This spectral measure is defined in a non-unique way if $A(t)$ is
in the limit circle case (see Section \ref{Section_spectr_repr}).
The moments of $d\rho^t(\lambda)$ are defined according to the
following rule.
\begin{equation*}
\label{Moment_eq_1} s_k(t)=\int_{-\infty}^\infty
\lambda^k\,d\rho^t(\lambda),\quad k=0,1,2,\ldots
\end{equation*}
In this section, we outline the derivation of the evolution of
$s_k(t)$ following \cite{MM12,MM13,MM17}.

We also consider the initial value problem for the finite Toda lattice associated with the $A^N(t)$:
\begin{equation}
\label{Toda_eq_fin}
\begin{cases}
\dot a_{N,\,n}(t)=a_{N,\,n}(t)\left(b_{N,\,n+1}(t)-b_{N,\,n}(t)\right),\\
\dot b_{N,\,n}(t)=2\left(a_{N,\,n}^2(t)-a_{N,\,n-1}^2(t)\right),
\end{cases},\quad
t\geqslant 0,\, n=1,2,\ldots,N.
\end{equation}
Here, we are looking for a solution that satisfies the initial condition
\begin{equation*}
\label{Toda_init_fin} a_{N,\,n}(0)=a_{N,\,n}^0,\quad
b_{N,\,n}(0)=b_{N,\,n}^0,\quad n=1,\ldots,N,
\end{equation*}
where $a_{N,\,n}^0,$ $b_{N,\,n}^0$ are real and $a_{N,\,n}^0>0,$
$n=1,2,\ldots,N$. It is a well known fact \cite{To,Te} that the
system (\ref{Toda_eq_fin}) is equivalent to (\ref{Toda_eqviv})
with the operator $A_N(t)$ given by $N\times N$ block of
(\ref{Jac_matr}) and $P_N(t)$ is defined by
\begin{align*}
&\left(P_N(t)f\right)_1=a_1(t)f_2,\\
&\left(P_N(t)f\right)_n=a_n(t)f_{n+1}-a_{n-1}(t)f_{n-1},\quad n=2,\ldots,N-1, \\
&\left(P_N(t)f\right)_N=a_{N-1}(t)f_{N-1}.
\end{align*}
The spectral measure $d\rho_N^t$  associated with $A^N(t)$ is
defined by (\ref{measure}) and has a form
\begin{equation*}
\label{Spec_mes} d\rho^t_N(\lambda)=\sum_{k=1}^N
\sigma_{N,\,k}^2(t)\delta\left(\lambda-\lambda_{N,\,k}(t)\right),
\end{equation*}
where values $\lambda_{N,\,k}(0)$, $\sigma_{N,\,k}(0)$,
$k=1,\ldots,N$ correspond to initial data $A^N(0)$. Moments of
this measure are denoted by $s_{N,k}(t)$, $k=0,1,\ldots$.

The authors adapt some results from \cite{M75} to find
$\lambda_{N,\,k}(t),\,\sigma_{N,\,k}(t)$:
\begin{proposition}
\label{Prop_eigen} The eigenvalues of the matrix $A^N(t)$ do not
depend on $t$:
\begin{equation*}
\lambda_{N,\,k}(t)=\lambda_{N,\,k}(0)=\lambda_{N,\,k},\quad
k=1,\ldots,N.
\end{equation*}
The norming coefficients are given by the \emph{Moser formula}:
\begin{equation*}
\label{Moser_f}
\sigma_{N,\,k}^2(t)=\frac{\sigma_{N,\,k}^2(0)e^{2\lambda_{N,\,k}t}}{\sum_{j=1}^N\sigma_{N,\,j}^2(0)e^{2\lambda_{N,\,j}t}},
\quad k=1,\ldots,N.
\end{equation*}
\end{proposition}
The proof is given in \cite[Section 2]{MM12}. Using Proposition
\ref{Prop_eigen} the authors obtain the recursive procedure for
determining the moments:
\begin{proposition}
The moments $s_{N,\,k}(t)$ of the measure $d\rho_N^t(\lambda)$
satisfy the following recurrent relation:
\begin{eqnarray}
\label{Main_syst} \dot
s_{N,\,k}(t)+\left(\ln{\left\{\|\Theta_N(t)\|^2\right\}}\right)'s_{N,\,k}(t)=2s_{N,\,k+1}(t),\quad
k=0,1,\ldots,\\
 \|
 \Theta_N(t)\|=\sqrt{\sum_{j=1}^N\sigma_{N,\,j}^2(0)e^{2\lambda_{N,\,j}t}}.\notag
\end{eqnarray}
\end{proposition}
The proof is given in \cite[Section 3]{MM12}. Since we know that
$s_{N, 0}(t) = 1$ for all $t$, we can use equation
(\ref{Main_syst}) to determine $s_{N,1}(t), s_{N,2}(t), \ldots,
s_{N,2N-2}(t)$ recursively. Then, we can use the fact that the set
of moments $s_{N,k}(t)$ for $k = 0, \ldots, 2N-2$ determines the
coefficients $a_{N,k}(t), b_{N,k}(t), a_{N,N}(t)$ for $k = 1,
\ldots, N-1$ by using the relationship with the response vector
given in equation (\ref{Resp_and_moments}) and the determinant
formulas (\ref{AK}) and (\ref{BK}).

The Jacobi matrix corresponding to initial data (\ref{Toda_init})
is denoted by $A_0$. We denote the $N\times N$ block of $A_0$ by
$A_0^{N}$.  The spectral measure of the operator corresponding to
$A_0^{N}$ is denoted by $d\,\rho_{N}^0(\lambda)$ and defined by
the formula (\ref{measure}). We know (see Section
\ref{Section_spectr_repr}) that $d\rho^0_N$ converges $*-$weakly
to $d\rho^0_{\alpha^*}$ as $N\to\infty,$ where
$d\rho^0_{\alpha^*}$ is a spectral measure of self-adjoint
extension $A^0_{\alpha^*}$ with $\alpha^*$ is defined by
(\ref{alpha_def}) in the limit circle case, or
$d\rho^0_{\alpha^*}$ is a spectral measure of $A_0$ in the limit
point case, see Section \ref{Section_spectr_repr}.

In \cite{MM17} the authors observed that if the support of
$d\,\rho_{\alpha^*}^0(\lambda)$ is semibounded, then it is
possible to pass to the limit as $N\to\infty$ in aforementioned
formulas.
\begin{proposition}
If the measure $d\rho_{\alpha^*}^0(\lambda)$ is such that
\begin{equation*}
\label{measure_restr}
\operatorname{supp}\left\{d\,\rho_{\alpha^*}^0(\lambda)\right\}\subset
(-\infty,M)
\end{equation*}
for some $M\in \mathbb{R}$, then there exist the limits
\begin{eqnarray}
s_k(t):=\lim_{N\to\infty}s_{N,\,k}(t)=\frac{\int_{\mathbb{R}}\lambda^ke^{2\lambda
t}\,d\rho_{\alpha^*}^0(\lambda)}{\int_{\mathbb{R}}e^{2\lambda
t}\,d\rho_{\alpha^*}^0(\lambda)},\quad k=0,1,\ldots,\label{SK_def}\\
\Omega_\alpha(t):=\lim_{N\to\infty}\|\Theta_N(t)\|^2=\int_{\mathbb{R}}e^{2\lambda
t}\,d\rho_{\alpha^*}^0(\lambda).\notag
\end{eqnarray}
Moreover, the moments $s_k(t)$ satisfy the recurrent relation
\begin{equation}
\label{Main_eq} \dot
s_{k}(t)+\left(\ln{\left\{\Omega_\alpha(t)\right\}}\right)'s_{k}(t)=2s_{k+1}(t),\quad
k=0,\ldots,\,\,t>0,
\end{equation}
where $s_0(t)=1$ for $t\geqslant 0.$
\end{proposition}
The proof is given in \cite[Proposition 7]{MM17}. It can be
demonstrated that the values $s_k(t)$ for any $t
> 0$, defined  using either equation (\ref{SK_def}) or equation
(\ref{Main_eq}), are actually the moments of a certain Borel
measure on the real line. This measure, in turn, corresponds to
the Jacobi matrix, which varies depending on the value of $t$.
Subsequently, the solution $a_n(t), b_n(t)$ for $n = 1, 2, \ldots$
to the Toda equations (\ref{Toda_eq}) with initial conditions
(\ref{Toda_init}) is determined by the moments $s_k(t), k = 0, 1,
\ldots$, using the relations (\ref{Resp_and_moments}) and the
determinant formulas (\ref{AK}) and (\ref{BK}).

\section{De Branges spaces associated with discrete dynamical
systems} \label{Section_DeBranges}

In \cite{MM3,MM15,MM16,MM21}, the authors established
relationships between  the Boundary control method \cite{BDAN87}
for one-dimensional systems and the de Branges method \cite{RR}.
In this section, we outline these relationships on the example of
the dynamic inverse problem for the system (\ref{Jacobi_dyn}). We
construct corresponding finite and infinite de Branges spaces. In
the first subsection, we provide information on de Branges spaces
according to \cite{RR}. The second section is based on \cite{MM3},
where finite-dimensional spaces (corresponding to blocks $A^N$ of
$A$) were constructed. In the third section, we demonstrate that
in the procedure of construction, the ``limit'' as $N\to\infty$ is
allowed if the operator $A$ (see Section
\ref{Section_spectr_repr}) is in the limit circle case. This
subsection is based on the work \cite{MM21} which was inspired by
the research of \cite{Berg1,Berg3,Y, Y2}, who studied properties
of operators associated with Hankel matrices.

\subsection{de Branges spaces.}

The entire function $E:\mathbb{C}\mapsto \mathbb{C}$ is called an
\emph{Hermite-Biehler function} if $|E(z)|>|E(\overline z)|$ for
$z\in \mathbb{C}_+$. We use the notation
$F^\#(z)=\overline{F(\overline{z})}$. The \emph{Hardy space} $H_2$
is defined by: $f\in H_2$ if $f$ is holomorphic in $\mathbb{C}^+$
and $\sup_{y>0}\int_{-\infty}^\infty|f(x+iy)|^2\,dx<\infty$. Then
the \emph{de Branges space} $B(E)$ consists of entire functions
such that
\begin{equation*}
B(E):=\left\{F:\mathbb{C}\mapsto \mathbb{C},\,F \text{ entire},
\,\frac{F}{E},\frac{F^\#}{E}\in H_2\right\}.
\end{equation*}
The space $B(E)$ with the scalar product
\begin{equation*}
[F,G]_{B(E)}=\frac{1}{\pi}\int_{\mathbb{R}}{\overline{F(\lambda)}}
{G(\lambda)}\frac{d\lambda}{|E(\lambda)|^2}
\end{equation*}
is a Hilbert space. For any $z\in \mathbb{C}$ the
\emph{reproducing kernel} is introduced by the relation
\begin{equation}
\label{repr_ker} J_z(\xi):=\frac{\overline{E(z)}E(\xi)-E(\overline
z)\overline{E(\overline \xi)}}{2i(\overline z-\xi)}.
\end{equation}
Then
\begin{equation*}
F(z)=[J_z,F]_{B(E)}=\frac{1}{\pi}\int_{\mathbb{R}}\overline{J_z(\lambda)}
{F(\lambda)}\frac{d\lambda}{|E(\lambda)|^2}.
\end{equation*}
We note that an Hermite-Biehler function $E(\lambda)$ determines
$J_z$ by (\ref{repr_ker}). The converse is also true
\cite{DMcK,DBr}:
\begin{theorem}
\label{TeorDB} Let $X$ be a Hilbert space of entire functions with
reproducing kernel $J_z(\xi)$ such that
\begin{itemize}

\item[1)] if $f\in X$ then $f^\#\in X$ and $\|f\|_X=\|f^\#\|_X$,

\item[2)] if $f\in X$ and $\omega\in \mathbb{C}$ such that
$f(\omega)=0$, then $\frac{z-\overline{\omega}}{z-\omega}f(z)\in
X$ and
$\left\|\frac{z-\overline{\omega}}{z-\omega}f(z)\right\|_{X}=\|f\|_{X}$,
\end{itemize}
then $X$ is a de Branges space based on the function
\begin{equation*}
E(z)=\sqrt{\pi}(1-iz)J_i(z)\|J_i\|_X^{-1}.
\end{equation*}
\end{theorem}

\subsection{Reachable sets and finite-dimensional de Branges spaces.}

According to \cite{AT} a spectral measure $d\rho_{t^*}(\lambda)$
corresponding to the operator $A$ (see Section
\ref{Section_spectr_repr}) give rise to the Fourier transform $F:
l^2\mapsto L_2(\mathbb{R},d\rho)$, defined by the rule:
\begin{equation*}
\label{JM_Four} (Fa)(\lambda)=\sum_{n=0}^\infty
a_k\phi_k(\lambda),\quad a=(a_0,a_1,\ldots)\in l^2,
\end{equation*}
where $\phi$ is a solution to (\ref{Phi_def}). The inverse
transform and Parseval identity have forms:
\begin{eqnarray}
a_k=\int\limits_{-\infty}^\infty
(Fa)(\lambda)\phi_k(\lambda)\,d\rho_{t^*}(\lambda),\notag\\
\sum_{k=0}^\infty \overline{a_k}{b_k}=\int\limits_{-\infty}^\infty
\overline{(Fa)(\lambda)}{(Fb)(\lambda)}\,d\rho_{t^*}(\lambda).\label{JM_parseval}
\end{eqnarray}
We assume that $T$ is fixed and $f\in \mathcal{F}^T$. Then, for
such a control and for $\lambda\in \mathbb{C}$, it is not hard to
see that we have the following representation for the Fourier
transform of the solution (\ref{Jac_sol_spectr_inf}) to the
equation (\ref{Jacobi_dyn}) at $t=T$ as follows:
\begin{equation*}
\left(Fu^f_{\cdot,T}\right)(\lambda)=\sum_{k=1}^T
\mathcal{T}_k(\lambda)f_{T-k},\quad \lambda\in \mathbb{C}.
\end{equation*}
Now we 
introduce the linear manifold of Fourier images of states of
dynamical system (\ref{Jacobi_dyn}) at time $t=T$, i.e. the
Fourier image of the reachable set:
\begin{equation*}
B_{A}^T\!\!:=\!F\mathcal{U}^T\!\!=\!\left\{\left(Fu^{J_Tf}_{\cdot,T}\right)(\lambda)\,|\,
J_T f\in \mathcal{F}^T\right\},
\end{equation*}
where the operator $J_T:\mathcal{F}^T\to \mathcal{F}^T$ is defined
as $(J_T f)_n=f_{T+1-n}$ for $n=1,\dots,T$.  It would be
preferable for us to use $C_T$ (see (\ref{C_T_eqn_ker}) for the
definition) instead of $C^T$.

We equip $B_{A}^T$ with the scalar product defined by the rule:
\begin{gather}
[F,G]_{B^T_{A}}:=\left(C_Tf,g\right)_{\mathcal{F}^t}, \quad F,G\in
B^T_{A},\label{JM_scalprod}\\
F(\lambda)=\sum_{k=1}^T
\mathcal{T}_k(\lambda)f_{k},\,G(\lambda)=\sum_{k=1}^T
\mathcal{T}_k(\lambda)g_{k},\,f,g\in \mathcal{F}^T.\notag
\end{gather}
Evaluating (\ref{JM_scalprod}) making use of (\ref{JM_parseval})
yields:
\begin{multline}
[F,G]_{B^T_{A}}=\left(C_Tf,g\right)_{\mathcal{F}^T}=\left(C^TJ_Tf,J_Tg\right)_{\mathcal{F}^T}
=\left(u^{J_Tf}_{\cdot,T},u^{J_Tg}_{\cdot,T}\right)_{\mathcal{H}^T}\label{DBr_Scal}\\
=\int\limits_{-\infty}^\infty
\overline{(Fu^{J_Nf}_{\cdot,N})(\lambda)}{(Fu^{J_Ng}_{\cdot,N})(\lambda)}\,d\rho_{t^*}(\lambda)
=\int\limits_{-\infty}^\infty
\overline{F(\lambda)}{G(\lambda)}\,d\rho_{t^*}(\lambda).
\end{multline}
Where the last equality is due to the finite speed of wave
propagation in (\ref{Jacobi_dyn}), (\ref{Jacobi_dyn_int}).

This shows that (\ref{DBr_Scal}) is a scalar product in
$B^T_{A}$. But we can say even more:
\begin{theorem}
The dynamical system with discrete time $(\ref{Jacobi_dyn})$ can
be used to construct the de Branges space
\begin{equation}
\label{De_Br_sp}
B_{A}^T:=\left\{\left(Fu^{J_Tf}_{\cdot,T}\right)(\lambda)\,|\, J_T
f\in \mathcal{F}^T\right\}=\left\{\sum_{k=1}^T
\mathcal{T}_k(\lambda)f_{k}\,|\, f\in \mathcal{F}^T\right\}.
\end{equation}
As a set of functions, it coincides with the space of Fourier
transforms of the states of the dynamical system
$(\ref{Jacobi_dyn})$ at time $t=T$ (or equivalently, the states of
$(\ref{Jacobi_dyn_int})$ with $N=T$ at the same time). In other
words, it is the set of Fourier transforms of reachable states,
and it consists of polynomials with real coefficients of degree at
most $N-1$. The norm in $B_{A}^T$ is defined using the connecting
operator.
\begin{equation*}
[F,G]_{B^T_{A}}:=\left(C_Tf,g\right)_{\mathcal{F}^T},\quad F,G\in
B^T_{A},
\end{equation*}
where
\begin{equation*}
F(\lambda)=\sum_{k=1}^T
\mathcal{T}_k(\lambda)f_{k},\,G(\lambda)=\sum_{k=1}^T
\mathcal{T}_k(\lambda)g_{k},\,f,g\in \mathcal{F}^T.
\end{equation*}
The reproducing kernel has a form
\begin{equation*}
J^T_z(\lambda)=\sum_{k=1}^T \mathcal{T}_k(\lambda)(j^z_T)_{k},
\end{equation*}
where $j^z_T$ is a solution to Krein-type equation
\begin{equation}
\label{Krein_eq}
C_Tj_T^z=(\mathcal{T}_{1}(z),
\mathcal{T}_{2}(z),\dots, \mathcal{T}_T(z))^*.
\end{equation}
\end{theorem}
Note (see Theorem \ref{Theor_Krein}) that control $J_T j_T^z$
drives the system (\ref{Jacobi_dyn}) to the special state at
$t=T$:
\begin{equation*}
\left(W^TJ_Tj_T^z\right)_i=\left(W^T_TJ_Tj_T^z\right)_i=\overline{\phi_i(z)},\quad
i=1,\ldots,T.
\end{equation*}
Thus, the reproducing kernel in the space $B^T_{A}$ is given by
\begin{multline*}
J_z^T(\lambda)=\sum_{k=1}^T \mathcal{T}_k(\lambda)(j^z_T)_{k}=\left(C_Tj^\lambda_T,j^z_T\right)_{\mathcal{F}^T}=\left(W^*_TW_Tj^\lambda_T,j^z_T\right)_{\mathcal{F}^T}\\
=\left(W_Tj_T^\lambda,W_Tj_T^z\right)_{\mathcal{H}^T}=\sum_{n=1}^T\overline{\phi_n(z)}{\phi_n(\lambda)}.
\end{multline*}

\begin{remark}
Due to the finite speed of wave propagation in (\ref{Jacobi_dyn}),
(\ref{Jacobi_dyn_int})  (cf. Remarks \ref{Rem_finite},
\ref{Rem11}), we can use the system (\ref{Jacobi_dyn_int}) with
$N=T$ and use the operator $C^T_T$ and the measure
$d\rho_T(\lambda)$ in formulas for the scalar product
(\ref{DBr_Scal}). We can also use any measure described in
Procedure \ref{Proc1}.
\end{remark}

\subsection{More on the connecting operator $C^T$ and infinite-dimensional de Branges spaces.}

Of course, it doesn't make sense to talk about what happens at
time $T=\infty$.  However, we can try to take the limit as $T$
approaches infinity in the formulas from the previous section. In
light of the relationships between the connecting operators $C^T$
and the Hankel matrices $S^T$ (as described in Proposition
\ref{CT_ST_prop}), using either $C^T$ or $S^T$ in the definition
of the norm in $B^T_A$ is equivalent. This is because we are
dealing with a finite-dimensional situation, and these matrices
differ by a factor. However, the properties of the
``infinite-dimensional'' versions of these operators differ, as
shown in references such as \cite{Berg1,Berg3,Y,Y2,MM21}.

We consider the  matrix $C$ formally defined by the product of two
matrices: $C=\left(W\right)^* W$, where
\begin{equation*}
W=\begin{pmatrix}
a_0 & w_{1,1} & w_{1,2} & \ldots \\
0 & a_0a_1 & w_{2,2} &  \ldots \\
\cdot & \cdot & \cdot & \cdot  \\
0 & \ldots & \prod_{j=1}^{k-1}a_j& \ldots \\
\cdot & \cdot & \cdot  & \cdot  \\
0 & 0 & 0  & \ldots
\end{pmatrix},
\end{equation*}
compare with (\ref{WT}) and
\begin{equation*}
C=
\begin{pmatrix}
r_0 & r_1& r_2 & \ldots & \ldots \\
r_1& r_0+r_2 & r_1+r_3 & \ldots &\ldots\\
r_2 & r_1+r_3 & \cdot & \cdot & \cdot \\
\ldots & \ldots & \ldots & \ldots & \ldots \\
\ldots & \ldots & \ldots & \ldots & \ldots
\end{pmatrix},\quad
\{C\}_{ij}=\sum\limits_{k=0}^{\max{\{i,j\}}-1}r_{|i-j|+2k},
\end{equation*}
compare with (\ref{C_T_repr}). We suggest that matrix $C $ should
play the same role in semi-infinite case as $C_T$ plays in
finite-dimensional situation.

The matrix $C_T$ is connected with the classical Hankel matrix
$S_T=J_TS^TJ_T$ by the following rule:
\begin{equation*}
C_T=\Lambda_T S_T\left(\Lambda_T\right)^*
\end{equation*}
where the matrix $\Lambda_T$ are defined in (\ref{LambdaMatr})

We rewrite the Krein equation (\ref{Krein_eq}) in terms of $S_T$:
counting (\ref{Resp_and_moments}) we get that
\begin{equation*}
S_T\Lambda_T^*j_T^\lambda=(1,\lambda,\dots,\lambda^{T-1})^*.
\end{equation*}

And introducing the vector $f_T^\lambda:=\Lambda_T^*j_T^\lambda$
we come to the equivalent form of (\ref{Krein_eq}):
\begin{equation}
\label{Krein_eq2} S_T f_T^\lambda=(1,\lambda,\dots,\lambda^{T-1})^*.
\end{equation}
Therefore the reproducing kernel in $B^T_{A}$ is  represented by
\begin{equation*}
J^T_z(\lambda)\!=\!\left(C_T
j^\lambda_T,j^z_T\right)_{\mathcal{F}^T}=
\left((1,\lambda,\dots,\lambda^{T-1})^*,\Lambda_T^*j^z_T\right)_{\!\!\!\mathcal{F}^T}\!\!=\!\sum_{n=0}^{T-1}f_k^z\lambda^k.
\end{equation*}

Krein equations in the form (\ref{Krein_eq}) and in the form
(\ref{Krein_eq2}) demonstrate the importance of the knowledge of
the invertability properties of operators $S_T$ and $C_T$ when $T$
goes to infinity. The authors in \cite{Berg1} studied this
question for $S_T$, below we answer the same question for $C_T$.

We introduce the notation:
\begin{equation*}
\beta_T:=\min\bigl\{l_k\,|\,l_k\text{ is an eigenvalue of
$C_T$},\, k=1,\ldots,T \bigr\}.
\end{equation*}

In what follows we mostly deal with $A$ in the limit-circle case.
The self-adjoint extensions and their spectral measures are
defined in the Section \ref{Section_spectr_repr}. By $dM(\lambda)$
we denote any of $d\rho_{\alpha}(\lambda)$ and by
$\left\{p_k(\lambda)\right\}_{k=0}^\infty$ denote the polynomials
introduced in Section \ref{Section_spectr_repr}.

\begin{theorem}
If the moment problem associated with sequence $\{s_k\}_{k=0}^\infty$ is
indetermined $($the matrix $A$ is in the limit circle case$)$ then
\begin{equation*}
\lim_{T\to\infty}\beta_T\geqslant
\left(\int\limits_{-1}^{1}m(x)\frac{dx}{\sqrt{1-x^2}}\right)^{-1},\quad
m(z)=\sum_{k=0}^\infty|p_k(z)|^2,
\end{equation*}
\end{theorem}
The proof of this theorem is given in \cite{MM21}.
\begin{remark}
We notice that the corresponding results for $S_{T}$ provide the
characterization of the limit point or limit circle case for
operator $A$. However, for $C_{T}$, the "if" part does not hold. A
straightforward example of this is the free Jacobi operator, where
$a_{k}=1$, $b_{k}=0$, for $k=1, 2, \ldots$. In this case, it's
easy to see that $C_{N}=I_{N}$ is an identity operator. As such,
for any $N$, $\beta_{N}=1$.
\end{remark}

Let $\beta^N$ denotes the maximal eigenvalue of $C_N$:
\begin{equation*}
\beta^N=\max\{l_k\,|\,l_k\text{ is eigenvalue of $C_N$},\,
k=1,\ldots,N \}.
\end{equation*}
\begin{lemma}
If $\beta^N\leqslant M$ then $A$ in the limit point case (the
moment problem associated with $\{s_k\}_{k=0}^\infty$ is
determinate).
\end{lemma}
The proof is given in \cite{MM21}.

The  matrix $C$ give rise to the (formally defined) operator
$\mathcal{C}$ in $l_2$:
\begin{equation*}
\left(\mathcal{C}f\right)_n:=\sum_{m=0}^\infty c_{nm}f_m,\quad
f\in l_2.
\end{equation*}
Then, without any a'priory assumptions on $c_{nm}$, only the
quadratic form
\begin{equation*}
C[f,f]:=\sum_{m,n\geqslant 0}c_{nm}\overline{f_m}{f_n}
\end{equation*}
is well-defined on the domain $D$ (\ref{domA}).

We always assume the positivity condition
\begin{equation*}
\sum_{m,n\geqslant 0}c_{nm}\overline{f_m}f_n\geqslant 0,\quad f\in
D,
\end{equation*}
which guarantees (see Theorem \ref{Th_char}) the existence of
Jacobi matrix $A$. The following result is valid

\begin{theorem}
The quadratic form $C[,]$ is closable in $l_2$ if one of the
following occurs:
\begin{itemize}
\item If the Jacobi operator $A$ is in the limit circle case
$($the moment problem associated with sequence $\{s_k\}$ is
indetermined$)$
\item If the Jacobi operator $A$ is bounded, i.e. for all
$k=1,2,\ldots$, $|a_k|,|b_k|\leqslant c$ for some $c\in
\mathbb{R}_+$ and its spectral measure is absolutely continuous
with respect to Lebesgue measure.
\end{itemize}
\end{theorem}
The proof is given in \cite{MM21}.

Now we assume that the sequence $\{s_k\}_{k=0}^\infty$  is
indeterminate (the matrix $A$ is in the limit circle case).  By
analogy with (\ref{De_Br_sp}) we introduce the linear manifold
\begin{equation*}
B_A^\infty:=\left\{\sum_{k=1}^\infty f_k \mathcal{T}_k(\lambda)\
\bigl{|}\ C[f,f]<\infty \right\}.
\end{equation*}
The scalar product in $B_A^\infty$ is given by the rule
\begin{eqnarray*}
[F,G]_{B_A^\infty}:=C[f,g]_{}=\int\limits_{-\infty}^\infty
\overline{F(\lambda)}{G(\lambda)}\,d\,M, \quad M\in
\mathcal{M}_H,\\
F,G\in B_A^\infty,\quad F(\lambda)=\sum_{k=1}^\infty f_k
\mathcal{T}_k(\lambda),\quad G(\lambda)=\sum_{k=1}^\infty g_k
\mathcal{T}_k(\lambda).
\end{eqnarray*}
The reproducing kernel is given by
\begin{equation*}
J^\infty_z(\lambda)=\sum_{n=1}^\infty
\overline{p_n(z)}{p_n(\lambda)}.
\end{equation*}
Then one can easily check the conditions of Theorem \ref{TeorDB}
to show that $B^\infty_A$ is a de Branges space.

\section{Weyl function for Jacobi operators}

In this section, we assume that the sequence of real numbers
$\{b_k\}_{k=0}^\infty$ and a sequence of positive numbers
$\{a_k\}_{k=0}^\infty$ with which we associate matrix $A$
(\ref{Jac_matr}) are both bounded:  $\sup_{n\geqslant
1}\{a_n,|b_n|\}<\infty$.

The operators $A^N$ and $A$ are defined in Sections \ref{Oper_fin}
and \ref{Section_spectr_repr}. It is important to note that, due
to conditions on coefficients, $A$ is in the limit point case
\cite{AH}.

The Weyl $m-$functions for $A$ and $A^N$ are defined by
\begin{eqnarray*}
m(\lambda)=\left( (A-\lambda)^{-1}e_1,e_1\right),\\
m^N(\lambda)=\left((A^N-\lambda)^{-1}e_1,e_1\right),
\end{eqnarray*}
$e_1=(1,0,0,\ldots)^T$.  By the spectral theorem, this definition
is equivalent to the following:
\begin{equation*}
m(\lambda)=\int_{\mathbb{R}}\frac{d\rho(z)}{z-\lambda},
\end{equation*}
where $d\rho=d\left(E_Ae_1,e_1\right)$ and $E_A$ is the
(projection-valued) spectral measure of the self-adjoint operator
$A$. This makes the Weyl function a crucial object in the study of
the  inverse spectral problems for one-dimensional systems, both
discrete and continuous \cite{AH,Te,BS1,AMR}.

In the work \cite{MMS}, a connection was found between the Weyl
functions $m(\lambda)$, $m^N(\lambda)$ and the response vectors
for systems (\ref{Jacobi_dyn}), (\ref{Jacobi_dyn_int}) defined in
(\ref{R_def}) and (\ref{R_def_int}). This connection provides a
new method for calculating Weyl functions and is also crucial in
the field of inverse spectral theory. The key idea is that there
is a relationship between special solutions of the dynamical
system (\ref{Jacobi_dyn}) and special solutions of the
corresponding ``spectral'' system, based upon a special
Fourier-type transform. For the case of the Schr\"odinger operator
on a half-line such connections were explored in
\cite{BS1,AMR,AM,MM2,MM3,MM22}.

To construct the appropriate Fourier transform, we consider two
solutions of the equation
\begin{equation}
\label{Difer_eqn0} \psi_{n-1}+\psi_{n+1}=\lambda \psi_n,
\end{equation}
satisfying the initial conditions:
\begin{equation*}
P_0=0,\,\,P_1=1,\quad Q_0=-1,\,\, Q_1=0.
\end{equation*}
Clearly, one has:
\begin{equation*}
Q_{n+1}(\lambda)=P_{n}(\lambda).
\end{equation*}
On introducing the notation $U_n(\lambda)=P_{n+1}(2\lambda)$, we
see that $U$ satisfies
\begin{equation*}
\label{Difer_eqnU} U_{n-1}(\lambda)+U_{n+1}(\lambda)=2\lambda
U_n(\lambda),\quad U_0=1,\,\,\ U_1=2\lambda.
\end{equation*}
So $U_n$ are Chebyshev polynomials of the second kind.

We are looking for the unique $l^2$-solution $S$ to equation
(\ref{Difer_eqn0}) which has the form:
\begin{equation*}
S_n(\lambda)=Q_n(\lambda)+m_0(\lambda)P_n(\lambda),
\end{equation*}
where $m_0$ is the Weyl function of $A_0$ (the special case of $A$
with  $a_n\equiv 1$, $b_n\equiv 0$). Consider the new variable
$z$, which is related to $\lambda$ by  the Joukowsky transform:
\begin{equation}
\label{Lam_z_conn} \lambda=z+\frac1z, \quad
z=\frac{\lambda-\sqrt{\lambda^2-4}}2=\frac{\lambda-i\sqrt{4-\lambda^2}}2,
\end{equation}
$z\in\mathbb D\cap\mathbb C_-$ for $\lambda\in\mathbb C_+$, where
$\mathbb{D}=\{z\,|\, |z|\leqslant 1\}$. Using the formula for
Chebyshev polynomials and the representation of the Weyl function
corresponding to finite matrix $A_0^N$ we can pass to the limit as
$N\rightarrow+\infty$ and obtain
$$
m_0(\lambda)=-z.
$$
Equation \eqref{Difer_eqn0} has two solutions, $z^n$ and $z^{-n}$.
Since $S\in l^2(\mathbb N)$ and $S_1(\lambda)=m_0(\lambda)=-z$, we
get
$$
S_n(\lambda)=-z^n.
$$
Let $\chi_{(a,b)}(\lambda)$ be the characteristic function of the
interval $(a,b)$. The spectral measure of the unperturbed operator
$A_0$ corresponding to (\ref{Difer_eqn0}) is
\begin{equation*}
d\rho(\lambda)=\chi_{(-2,2)}(\lambda)\frac{\sqrt{4-\lambda^2}}{2\pi}\,d\lambda.
\end{equation*}
Then, the Fourier transform $F: l^2(\mathbb{N})\mapsto
L_2((-2,2),\rho)$ is defined as follows: for
$(v_n)\in l^2(\mathbb{N})$:
\begin{equation*}
F(v)(\lambda)=V(\lambda)=\sum_{n=0}^\infty S_n(\lambda)v_n.
\end{equation*}
The inverse transform is given by:
\begin{equation*}
v_n=\int\limits_{-2}^2V(\lambda)S_n(\lambda)\,d\rho(\lambda).
\end{equation*}

We recall that finite and semi-infinite Jacobi operators are
introduced in Sections \ref{Oper_fin},
(\ref{Section_spectr_repr}). The main result of \cite{MMS} is the
following.
\begin{theorem}
If the coefficients in the semi-infinite Jacobi
matrix operator $A$ or the finite Jacobi matrix operator $A^N$
satisfy the condition $\sup_{n \geqslant 1} \{a_n, |b_n|\}
\leqslant B$, then the Weyl functions $m$ and $m^N$ can be
represented as
\begin{eqnarray*}
m(\lambda)=-\sum_{t=0}^\infty z^tr_t,\label{M_inf}\\
m^N(\lambda)=-\sum_{t=0}^\infty z^tr^N_t\label{M_fin}.
\end{eqnarray*}
Here, $(r_0, r_1, r_2, \ldots)$ and $(r_0^N, r_1^N, r_2^N,
\ldots)$  represent response vectors, which are the kernels of
dynamic response operators related to (\ref{Jacobi_dyn}) and
(\ref{Jacobi_dyn_int}).  Variables $\lambda$ and $z$ are related
by Joukowsky transform (\ref{Lam_z_conn}). These representations
hold for $\lambda\in D$, where $D\subset \mathbb{C}_+$  is defined
as follows. Let $R:=3B+1$, then
\begin{eqnarray*}
\label{D_reg} D:=\left\{x+iy\,\Big|\, x>
\left(R+\frac{1}{R}\right)\cos{\phi},\right.\\
\left.y>0,\,
y>\left(\frac{1}{R}-R\right)\sin{\phi},\,\phi\in(\pi,2\pi)\right\}.\notag
\end{eqnarray*}
\end{theorem}
The proof of this theorem is given in \cite[Theorem 1]{MMS}.

\section{Discrete dynamical systems with continuous time}
\label{Section_KreinString}

In this section, we solve the dynamic IP for a system defined by
equation (\ref{DnSst}). Here, the ``space'' variable is
represented by discrete values, while the ``time'' variable is
continuous.

In the first part of this section, we study the forward problem
for equation (\ref{DnSst}) and derive representations for the
operators used in the BC method. The authors in \cite{MM14}
proposed two methods for recovering $A^N$. The first method relies
on trace-type formulas, while the second method exploits the
possibility of extracting spectral data from dynamic data related
to $A^N$. The latter method allows us to reduce the problem to a
well-known one, as described in Remark \ref{Rem2} and Remark
\ref{Rem3}. In this review, we primarily focus on the method based
on Krein equations. Additionally, we provide a characterization of
inverse data. In the second part of this section, we introduce a
dynamic IP for a Krein string and delve deeper into a particularly
significant case: a finite Krein-Stieltjes string. Finally, in the
last part of this section, we investigate the problem of
approximating a finite string with constant density by a
weightless string where the masses are evenly distributed.

\subsection{Inverse dynamic problem for finite Jacobi matrix.}

We fix some $T>0$ and introduce the \emph{outer space} of the
system (\ref{DnSst}), the space of controls:
$\mathcal{F}^T:=L_2(0,T)$ with the scalar product $f,g\in
\mathcal{F}^T,$
$\left(f,g\right)_{\mathcal{F}^T}=\int_0^Tf(t)g(t)\,dt$. The
solution of (\ref{DnSst}) with such a control we denote by $u^f$.

We recall that (\ref{Phi_def}) determines the set of polynomials
$\{1,\phi_2(\lambda),\ldots,\phi_N(\lambda),\phi_{N+1}(\lambda)\}$,
where $\phi_k(\lambda)$ is a polynomial of degree $k-1$. By
$\lambda_1,\ldots,\lambda_{N}$ we denote the roots of the equation
$\phi_{N+1}(\lambda)=0$. It is known \cite{AH} that they are real
and distinct. The vectors
$\phi(\lambda),\,\left\{\phi^k\right\}_{k=1}^N$ and weights
$\left\{\omega_k\right\}_{k=1}^N$ are introduced in (\ref{Phi_k}),
(\ref{Ortog}).
Thus, $\phi^k$ are non-normalized eigenvectors of $A^N$,
corresponding to eigenvalues $\lambda_k$:
\begin{equation*}
A^N\phi^k=\lambda_k\phi^k,\quad \quad k=1,\ldots,N.
\end{equation*}
We call by \emph{spectral data} and spectral function $\rho$ the
following objects:
\begin{equation*}
\label{measure_JM}
\left\{\lambda_i,\omega_i\right\}_{i=1}^{N},\quad
\rho(\lambda)=\sum_{\{k:\lambda_k<\lambda\}}\frac{1}{\omega_k}.
\end{equation*}
Note the property:
\begin{equation}
\label{rho_pr} \sum_{k=1}^N\frac{1}{\omega_k}=1.
\end{equation}

The standard application of a Fourier method yields the
representation formula:
\begin{lemma}
The solution to (\ref{DnSst}) admits the spectral representation
\begin{equation}
\label{Sol_spec_repr_JM} u^f(t)=\sum_{k=1}^{N}h_k(t)\phi^k,\quad
u^f(t)=\int_{-\infty}^\infty\int_0^t
S(t-\tau,\lambda)f(\tau)\,d\tau\phi(\lambda)\,d\rho(\lambda).
\end{equation}
where
\begin{eqnarray*}
h_k(t)=\frac{1}{\omega_k}\int_0^t
f(\tau)S_k(t-\tau)\,d\tau,\notag\\
S(t,\lambda)=\begin{cases}
\frac{\sin{\sqrt{\lambda}t}}{\sqrt{\lambda}},\quad
\lambda>0,\\
\frac{\operatorname{sh}{\sqrt{|\lambda|}t}}{\sqrt{|\lambda|}},\quad
\lambda<0,\\
t,\quad \lambda=0,
\end{cases},\quad S_k(t)=S(t,\lambda_k).\label{S_k}
\end{eqnarray*}
\end{lemma}

The \emph{response operator} $R^T: \mathcal{F}^T\mapsto
\mathcal{F}^T$ is introduced by formula (\ref{RespJM}). Making use
of (\ref{Sol_spec_repr_JM}) implies the representation for $R^T$:
\begin{equation*}
\label{Resp_def}
\left(R^Tf\right)(t)=u_1^f(t)=\sum_{k=1}^{N}h_k(t)=\int_0^t
r(t-s)f(s)\,ds,
\end{equation*}
where
\begin{equation}
\label{resp_func_JM}
r(t)=\sum_{k=1}^{N}\frac{1}{\omega_k}S_k(t),
\end{equation}
is called a \emph{response function}. Note that the operator $R^T$
is a natural analog of a dynamic Dirichlet-to-Neumann operator
\cite{B07,B17} in continuous, and (\ref{resp_oper}) in discrete
cases.

A key difference between system (\ref{DnSst}) and the ones
described in (\ref{Jacobi_dyn}) and (\ref{Jacobi_dyn_int}) is an
infinite wave propagation speed in (\ref{DnSst}):
\begin{remark}
Formula  (\ref{Sol_spec_repr_JM}) implies  an infinite wave
propagation speed in (\ref{DnSst}) the latter means that for the
wave initiated by the control $F=(f,0,\ldots,0)$ in (\ref{DnSst})
we have that $u^f_n(\tau)\not=0$ for any $n=1,2,\ldots,N$ and any
$\tau>0$, without special assumptions on $f$.
\end{remark}
(Compare this with Remark \ref{Rem_finite}). The infiniteness
speed of wave propagation in (\ref{DnSst}) significantly
complicates the process of solving the dynamic IP.

The natural setup for dynamic IP involves finding the matrix $A^N$
based on the knowledge of the response operator $R^T$ for
arbitrarily fixed $T>0$. This is equivalent by using the
information about $r(t)$, which is given by (\ref{resp_func_JM})
on the time interval $(0,T)$. The possibility of choosing any $T$
greater than zero is due to the infinite speed of wave propagation
in the system (\ref{DnSst}). This fact follows directly from the
representation (\ref{Sol_spec_repr_JM}).

The \emph{inner space}, i.e. the space of states of (\ref{DnSst})
is denoted by $\mathcal{H}^N:=\mathbb{R}^{N}$, indeed for any
$T=N>0$ we have that $u^f(T)\in \mathcal{H}^N$. The metric in
$\mathcal{H}^N$ is given by
\begin{equation*}
(a,b)_{\mathcal{H}^N}=\left(a,b\right)_{\mathbb{R}^{N}},\quad
a,b\in\mathcal{H}^N.
\end{equation*}
The \emph{control operator} $W^T: \mathcal{F}^T\mapsto
\mathcal{H}^N$ is introduced by the rule:
\begin{equation*}
W^Tf=u^f(T).
\end{equation*}
Due to (\ref{Sol_spec_repr_JM}) we have that
$W^Tf=\sum_{k=1}^{N}h_k(T)\phi^k$. We set up the following
\emph{control problem}:  for a fixed state $a\in \mathcal{H}^N$,
$a=\sum_{k=1}^{N}a_k\phi^k$, we seek a control $f\in\mathcal{F}^T$
driving the system (\ref{DnSst}) to a prescribed state:
\begin{equation}
\label{control_JM} W^Tf=a.
\end{equation}
Using (\ref{Sol_spec_repr_JM}) we see that the equality
(\ref{control_JM}) is equivalent to the following moment problem:
to find $f\in \mathcal{F}^T$ such that
\begin{equation*}
a_k=h_k(T)=\frac{1}{\omega_k}\int_0^T
f(\tau)S_k(T-\tau)d\tau,\quad k=1,\ldots,N.
\end{equation*}
Clearly, a solution to this moment problem exists, but is not
unique: it is a simple consequence of the fact that the outer
space is infinite-dimensional, and the space of states is
finite-dimensional. We introduce the subspace
\begin{equation*}
\mathcal{F}^T_1={\operatorname{span}\left\{S_k(T-t)\right\}_{k=1}^N}.
\end{equation*}
Then the following result on the boundary controllability of
(\ref{DnSst}) is valid:
\begin{lemma}
\label{LemmaCont_JM} The operator ${W}^T$ maps $\mathcal{F}^T_1$
onto $\mathcal{H}^N$ isomorphically.
\end{lemma}
The proof of this lemma is given in \cite[Lemma 2.2]{MM14}.

The \emph{connecting operator} $C^T:\mathcal{F}^T\mapsto
\mathcal{F}^T$ is defined by the rule
$C^T:=\left(W^T\right)^*W^T$, so by the definition for $f,g\in
\mathcal{F}^T$ one has
\begin{equation*}
\label{C_T_JM}
\left(C^Tf,g\right)_{\mathcal{F}^T}=\left(u^f(T),u^g(T)\right)_{\mathcal{H}^N}=\left(W^Tf,W^Tg\right)_{\mathcal{H}^N}.
\end{equation*}
It is crucial in the BC method that $C^T$ can be expressed in
terms of inverse data:
\begin{theorem}
The connecting operator admits the representation in terms of
dynamic inverse data:
\begin{equation}
\label{CT_dyn_repr_JM}
\left(C^Tf\right)(t)=\frac{1}{2}\int_0^T\int_{|t-s|}^{2T-s-t}r(\tau)\,d\tau
f(s)\,ds,
\end{equation}
and in terms of spectral inverse data:
\begin{equation}
\label{CT_sp_repr_JM} \left(C^Tf\right)(t)=\int_0^T
\sum_{k=1}^{N}\frac{1}{\omega_k}S_k(T-t)S_k(T-s)f(s)\,ds.
\end{equation}
\end{theorem}
Proof is given in \cite[Theorem 2.3]{MM14}.

\begin{remark}
The formula (\ref{CT_sp_repr_JM}) implies
that $\mathcal{F}^T_1=C^T\mathcal{F}^T$, that is the space
$\mathcal{F}^T_1$ is completely determined by inverse data. The
important properties of functions from $\mathcal{F}^T_1$ is that
\begin{equation*}
\label{FT_proper_JM} f(T)=0,\quad \text{for all}\,\, f\in
\mathcal{F}^T_1.
\end{equation*}
\end{remark}

By $f_k^T\in \mathcal{F}^T_1$ we denote controls, driving the
system (\ref{DnSst}) to prescribed \emph{special states}
\begin{equation*}
d_k\in \mathcal{H}^N, \,d_k=\left(0,\ldots,1,\ldots,0\right),\quad
k=1,\ldots,N.
\end{equation*}
In other words, $f^T_k$ are the solutions to the equations
$W^Tf_k^T=d_k$, $k=1\,\ldots,N$. By Lemma \ref{LemmaCont_JM}, we
know that these controls exist and are unique. It is important that
they can be found as the solutions to the Krein equations.

\begin{theorem} The control $f_1^T$ can be found as the solution to
the following equation:
\begin{equation*}
\label{Cont_f1_JM} \left(C^Tf_1^T\right)(t)=r(T-t),\quad 0<t<T.
\end{equation*}
The controls $f^T_k$ satisfy the system:
\begin{equation}
\label{Cont_system_JM}
\begin{cases}
-\left(C^Tf_1^T\right)''=b_1C^Tf^T_1+a_1C^Tf^T_{2},\\
-\left(C^Tf_k^T\right)''=a_{k-1}C^Tf^T_{k-1}+b_kC^Tf^T_k+a_kC^Tf^T_{k+1},\quad
k=2,\ldots,N-1,\\
-\left(C^Tf_N^T\right)''=a_{N-1}C^Tf^T_{N-1}+b_NC^Tf^T_N,
\end{cases}
\end{equation}
\end{theorem}
Proof is given in \cite[Theorem 2.5]{MM14}.

To recover unknown coefficient of $A$ we propose the following
procedure. First we note that by the definition of controls
$f^T_k$:
\begin{equation*}
\label{F_orthog}
\left(C^Tf^T_i,f^T_j\right)_{\mathcal{F}^T}=\delta_{ij}.
\end{equation*}
Thus, multiplying the first line in (\ref{Cont_system_JM}) by
$f^T_1$ in $\mathcal{F}^T$, we obtain that
\begin{equation*}
b_1=-\left(\left(C^Tf^T_1\right)'',f^T_1\right)_{\mathcal{F}^T}.
\end{equation*}
Then from the first line in (\ref{Cont_system_JM}) we have that
\begin{equation*}
a_1C^Tf^T_{2}=-\left(C^Tf_1^T\right)''-b_1C^Tf^T_1.
\end{equation*}
Inverting $C^T$ in $\mathcal{F}^T_1$, we can find the quantity
$a_1f^T_{2}$. Then evaluating the quadratic form
\begin{equation*}
\left(C^Ta_1f^T_2,a_1f^T_2\right)_{\mathcal{F}^T}=a_1^2\left(C^Tf^T_2,f^T_2\right)_{\mathcal{F}^T}=a_1^2,
\end{equation*}
we recover $a_1$ and $f^T_2$. Continuing this procedure we obtain
that
\begin{eqnarray*}
b_n=-\left(\left(C^Tf^T_n\right)'',f^T_n\right)_{\mathcal{F}^T},\\
a_n=-\left(\left(C^Tf^T_n\right)'',f^T_{n+1}\right)_{\mathcal{F}^T}=-\left(\left(C^Tf^T_{n+1}\right)'',f^T_n\right)_{\mathcal{F}^T}.
\end{eqnarray*}

Now we can formulate the result on the characterization of inverse
data:
\begin{theorem}
The function $r$ is a response function for the dynamical system
(\ref{DnSst}) if and only if it has a form (\ref{resp_func_JM})
where coefficients $\omega_k,$ $k=1,\ldots,N$ satisfy
(\ref{rho_pr}), operator $C^T$ defined by (\ref{CT_dyn_repr_JM})
is a finite-rank operator, $C^T\mathcal{F}^T:=\mathcal{F}^T_1$,
and $C^T|_{\mathcal{F}^T_1}$ is an isomorphism.
\end{theorem}
The proof is given in \cite[Theorem 2.9]{MM14}.

The corresponding de Branges space for the system (\ref{DnSst}) is
constructed in \cite{MM16} using the method from Section
\ref{Section_DeBranges}.

\subsection{Special case: Krein-Stieltjes string}

In this section, we describe the setting of the dynamic IP for a
Krein string, which is a historically very important object
\cite{Kr1,Kr2,KK}. We set up the IP and consider in more detail
the special case of a finite Krein-Stieltjes string.

Let's consider a non-decreasing bounded function $M(x)$ defined on
the interval $(0, l)$. Here, $M(x)$ represents the mass of the
string segment $(0, x)$. If the function $M$ is differentiable,
then we can define the string density $\rho(x)$ as $\rho(x) =
M'(x)$. Following \cite{DMcK, KK}, we introduce a domain that
consists of \emph{continued functions}
\begin{multline*}
    D_M:=\Bigl\{\left[u(x),u'_-(0),u'_+(l)\right]\,\bigl|\quad
    u(x)=a+bx-\int_0^x(x-s)g(s)\,dM(s);\\
    u'_-(0)=b;\quad u'_+(l)=b-\int_0^lg(s)\,dM(s) \Bigr\},
\end{multline*}
 where $ u '_- (0), u' _ + (l) $ are left and right derivatives, $ g $ is
 $ M $-summable function. Define the \emph{generalized differential Krein operation}
 $ l_M [u] $ on $ D_M $ by the rule
\begin{equation*}
    \label{L_oper} l_M[u]:=g(x),\quad \text{$M-$almost everywhere}.
\end{equation*}
   Note that if the function $ M $ is differentiable and $ M '(x)> 0 $, then
    $l_m[y] = -\frac{1}{M'(x)} \frac{d^2 y(x)}{dx^2}$.

We fix $ T> 0 $ and consider the following initial-boundary value
problem:
\begin{equation}
    \label{L_eq}
    \begin{cases}
        u_{tt}(x,t)+l_M[u]=0,\quad 0<x<l,\, 0<t<T,\\
        u(x,0)=0,\, u_t(x,0)=0,\quad 0\leqslant x\leqslant l,\\
        u(0,t)=f(t),\,  u(x,l)=0\quad 0\leqslant t\leqslant T,
    \end{cases}
\end{equation}
where $ f \in L_2 (0, T) $  is a \emph{boundary control}.
The dynamic inverse data is given by the \emph{response operator}
$R^T: L_2(0,T)\mapsto L_2(0,T)$ defined by
\begin{equation*}
\left(R^Tf\right)(t):=u^f_x(0,t),\quad t\in (0,T).
\end{equation*}
When $M'(x) = \rho(x) \ge \delta > 0$ for $0 \le x \le l$, and
$\rho \in C^1$,  then the system (\ref{L_eq}) corresponds to an
initial-boundary value problem for the wave equation. The dynamic
IP for such a system were studied in \cite{B08,BL}. However, the
case of a general density has not been widely discussed in the
literature. The reason for this is that the wave propagation speed
in systems with a ``bad'' density is infinite, making the usual
dynamic methods, such as the BC method, inapplicable in their
current form. Therefore, these methods need to be modified to
account for this.

\begin{lemma}
The initial boundary value problem (\ref{L_eq}) is equivalent to
the following integral equation:
\begin{multline*}
\int_0^x
z(x-z)u(z,t)\,dM(z)=x\int_0^t(t-s)\left[f(s)+u(x,s)\right]\,ds\\
-2\int_0^t(t-s)\int_0^x u(y,s)\,dy\,ds,\quad u(x,l)=0.
\end{multline*}
\end{lemma}
The proof of this Lemma is given in \cite[Section 2]{MM11}.


We assume that the string is a Krein-Stieltjes string, i.e. the
mass $ M $ is a piecewise constant function. Let
$0=x_0<x_1<x_2<\ldots<x_{N-1}<x_N=l$, $m_i>0,$ $i=1,\ldots,N-1$,
$l_i=x_i-x_{i-1}$, $i=1,\ldots N,$ and the density $ dM $ has the
form $dM(x)=\sum_{i=1}^{N-1}m_i\delta(x-x_i)$. In this case, we
find that
\begin{align*}
&l_M[u](x)=\frac{u'(x+0)-u'(x-0)}{m_x},\\
&m_x:=\begin{cases} M(+0)-M(0),\quad x=0,\\
M(x+0)-M(x-0),\quad 0<x<l,\\
M(l)-M(l-0),\quad x=l,
\end{cases}
\end{align*}
and, therefore, $ u '(x) = 0 $, $ x \in (x_ {i-1}, x_i) $ and
introducing the notation $ u_i (t) = u (x_i, t) $, we see that the
initial-boundary value problem (\ref{L_eq}) is equivalent to the
initial-boundary value problem for a vector function (we keep the
same notation) $u(t)=\left(u_1(t),\ldots,u_{N-1}(t)\right)$:
\begin{equation}
\label{DynSyst}
\begin{cases}
Mu_{tt}=-Au+\widetilde{f},\quad t> 0,\\
u(0)=0,\, u_t(0)=0,
\end{cases}
\end{equation}
where
\begin{eqnarray*}
A=\begin{pmatrix} b_1 & a_1&0&   \ldots & 0\\
a_1& b_2 & a_2& \ldots & 0\\
\ldots & \ldots& \ldots & \ldots \\
0 & 0& \ldots&  b_{N-2} &a_{N-2}\\
0 & 0& \ldots& a_{N-2} &b_{N-1}
\end{pmatrix},\\
M=\begin{pmatrix} m_1 & 0&  \ldots & 0\\
0& m_2 &  \ldots & 0\\
\ldots & \ldots& \ldots & \ldots \\
0 & \ldots& 0 & m_{N-1}
\end{pmatrix}, \,
\widetilde f=\begin{pmatrix} \frac{f}{l_1}\\0\\
\ldots\\0\end{pmatrix},
\end{eqnarray*}
and the elements of the matrix $ A $ are defined as
\begin{equation*}
\label{coeff} a_i=\frac{1}{l_{i+1}},\quad
b_i=-\frac{l_i+l_{i+1}}{l_il_{i+1}}.
\end{equation*}
We note that in this case $A<0$.

The solution of the initial-boundary value problem (\ref{L_eq}) is
defined as a continuous function according to the following rule:
\begin{equation}\label{con_v_f}
u^f(x,t)=\begin{cases}  f(t),\quad x=0,\\
u_i(t), \quad x=x_i,\quad i=1,\ldots,N-1,\\
0,\quad x=l,\\
\text{affine function},\quad x\in (x_{i-1},x_i).
\end{cases}
\end{equation}
The solution to (\ref{DynSyst}) is denoted by $ u^ f (t) $. In
what follows we  fix  $T>0$.

\emph{Dynamic inverse problem} consists in reconstructing the
string, that is, the coefficients $ l_N $, $ m_i $, $ l_i $, $ i =
1, \ldots, N-1 $ or, equivalently, the elements of the matrices $
A $ and $ M $, based on the \emph{response operator} $ R^ T $
defined by the rule
\begin{equation*}
\label{Resp} \left(R^Tf\right)(t):=u^f_1(t),\quad 0<t<T.
\end{equation*}
The IP for the system (\ref{DynSyst}) is solved in
\cite{MM14,MM11} using methods from the previous section.

\begin{remark}
We note that the BC method \cite{BDAN87,B97,B07} was developed
based on a local approach to solving the dynamic IP for an
inhomogeneous string \cite{BL}. That is why it is crucial to find
a solution for the dynamic IP for a string with arbitrary mass
distribution $M$.
\end{remark}

\subsection{Approximation of constant density by point masses}

In this section, we consider a specific scenario in which a string
represents the interval $(0, 1)$ with a constant density $\rho(x)
= M'(x) = 1$. We then consider an approximation of this string by
a weightless string with uniformly distributed point masses. For
this specific case, the initial-boundary value problem
(\ref{L_eq}) takes the following form:
\begin{equation}
\label{eq_unit}
\begin{cases}
v_{tt}(x,t)-v_{xx}=0,\quad 0<x<1,\, 0<t<T,\\
v(x,0)=0,\, v_t(x,0)=0,\quad 0\leqslant x\leqslant 1,\\
v(0,t)=f(t),\,  u(x,1)=0\quad 0\leqslant t\leqslant T.
\end{cases}
\end{equation}

Let us emphasize that in this system the speed of wave propagation
is finite and equal to one. As an approximation, we use the
identical point masses and lengths: $ l_i = 1 / N $, $ m_i = 1 / N
$, $ i = 1, \ldots, N-1 $ and the wave propagation velocity in an
approximate system (\ref{DynSyst}) is infinite. Our goal is to
obtain the effect of the appearance of a finite wave propagation
velocity at $ N \to \infty. $

In this situation, the matrices $A$ and $M$ in (\ref{DynSyst}) have the following forms:
\begin{equation}
A=N\begin{pmatrix} -2 & 1&   \ldots & 0\\
1& -2 &  \ldots & 0\\
\ldots & \ldots& \ldots &1 \\
0 &  \ldots & 1 &-2
\end{pmatrix},\quad 
M=\frac{1}{N}\begin{pmatrix} 1 & 0&  \ldots & 0\\
0& 1 &  \ldots & 0\\
\ldots & \ldots& \ldots & \ldots \\
0 & \ldots& 0 & 1
\end{pmatrix}.\label{MN}
\end{equation}

In \cite{MM11} the authors derived the formula for the solution of
(\ref{DynSyst}) with specified matrices.

\begin{proposition}\label{p2}
The response function $ r_N (t) $ for a dynamic system
(\ref{DynSyst})  with matrices defined as (\ref{MN}) (for the
Krein-Stieltjes strings with uniformly distributed identical point
masses) converges in the sense of distributions to
$$
r_N(t)\to  \delta(t) 
\quad\text{ as }N\to\infty.
$$
\end{proposition}
Where the convergence is understood as in \cite[Chapter 2]{V}, the
proof of this proposition is given in \cite[Proposition 2]{MM11}.

The Krein-Stieltjes string with uniformly distributed identical
point  masses serves as  a discrete approximation of a unit
density string. For the initial-boundary value problem for a
string with unit density (\ref{eq_unit}), it is known that the
response function takes the form $ r (t) = - \delta'(t) $.
According to Proposition \ref{p2},  the response function of the
approximate system does not converge to that of the original
system. This is due to the definition of the response function
(the response operator) in an approximate (discrete) system. In
the original system (\ref{eq_unit}), the response operator is
defined as $R^Tf(t) = v_x^f(0, t)$, which represents the value of
the derivative of the solution at the point $x = 0$. However, in
the approximate system  the reaction operator takes the form
$R^Tf(t) = u_1^f(t)$, which represents the value of the solution
in the first channel at the point $x_1 = 1/N$.

We can introduce another definition of \emph{corrected response}:
namely we set
\begin{equation}\label{pr}
\tilde r_N=\frac{u_1^f(t)-u_0^f(t)}{1/N},
\end{equation}
 then we can prove the following
\begin{proposition}\label{p3}
The corrected response function $ \tilde r_N (t) $ (\ref{pr}) for
the dynamical system (\ref{DynSyst}) with matrices  (\ref{MN})
(for the Krein-Stieltjes string with uniformly distributed
identical point masses) converges in the sense of distributions to
the response function for the system (\ref{eq_unit})
$$
\tilde r_N(t)\to  r(t)=-\delta'(t), \quad\text{ if }N\to\infty.
$$
\end{proposition}

\begin{proposition}\label{p4}
The solution $u^N$ (\ref{con_v_f}) of the initial-boundary value
problem (\ref{L_eq}) with matrices (\ref{MN}) converges in the
sense of distributions
$$
u^N(x,t)\to v^\delta(x,t)=\delta(x-t)\quad\text{ if
}N\to\infty,
$$
where $v^\delta(x,t)$\,--\, the solution of the problem
(\ref{eq_unit}).
\end{proposition}
The proofs of Propositions \ref{p3}, \ref{p4} are given in
\cite[Propositions 3,4]{MM11}.

The results mentioned above demonstrate that a "simple" system
(\ref{DynSyst}) can potentially be used to approximate "complex"
system (\ref{L_eq}). Similar approaches have been employed in
\cite{DBr, RR} to approximate complex systems with simpler
finite-dimensional systems. These simpler systems are easier to
analyze, making them valuable tools for understanding and managing
more complicated systems.

\section{Application to numerical simulation}

In this section we consider the idea of using discrete systems as
a way to approximate continuous ones.

First, we  derive the discrete version of the system
(\ref{Wave_eqn})-(\ref{Init}), that simulate the wave propagation
in the continuous model. It is important to note that trying to
create this discrete version by simply replacing the derivatives
with the corresponding differences does not work.

To understand where the equations (\ref{Wave_eqn})-(\ref{Init})
come from, let us first introduce the concepts of kinetic and
potential energies:
\begin{eqnarray*}
T(t)=\frac{1}{2}\int_\Omega u_t^2(x,t)\,dx, \label{Kin_en} \\
U(t)=\frac{1}{2}\int_\Omega u_x^2(x,t)\,dx. \label{Pot_en}
\end{eqnarray*}
According to the principles of classical mechanics, known as
Hamilton's  principle of least action, a system moves from state 1
at time $t_1$ to state 2 at time $t_2$ in such a way that the
variation of the action functional vanishes
\begin{equation*}
S[u]= \int_{t_1}^{t_2} L(t)\,dt,\quad L=T-U\text{ (Lagrangian)},
\end{equation*}
 along the actual path of movement. It means that
$$
\delta S[u,h]=0 \quad \text{for all}\quad h: h|_{t=t_1}=0,
h|_{t=t_2}=0, h|_\Gamma=0.
$$
Calculating the variation of the  functional and using integration
by parts, we arrive at the system (\ref{Wave_eqn})--(\ref{Init}).
A similar approach was proposed in \cite{AChMM} for discrete
systems, where a problem of  boundary controllability for the
discrete system was also considered.

In the last section, it is suggested to replace the system
(\ref{Parab_eqn}) with a discrete one related to the matrix $A^N$:
\begin{equation}
\label{Heat}
\begin{cases}
v_{n,t+1}-a_nv_{n+1,t}-a_{n-1}v_{n-1,t}-b_nv_{n,t}=0,\,\, t\in \mathbb{N}\cup\{0\},\,\, n\in \mathbb{N} \\
v_{n,\,0}=0,\quad n\in \mathbb{N}, \\
v_{0,\,t}=f_t,\quad 
\quad t\in \mathbb{N}\cup\{0\},
\end{cases}
\end{equation}
where $f=(f_0,f_1,\ldots)$\,---\,\emph{boundary control}, $a_0=1$.
Solution to (\ref{Heat}) is defined by $v^f$. Note that, compared
to (\ref{Jacobi_dyn}),  we replace the term
$u_{n,\,t+1}-u_{n,\,t11}$, which corresponds to the second
derivative $u_{tt}$, with the term $v_{n,t+1}$, which corresponds
to the first derivative $v_t$. We provide the details of the
solution of IP for (\ref{Heat}) following \cite{MM20,MMSh}.

\subsection{Discrete dynamical systems on graphs.}

In this section, instead of using a continuous metric graph
$\Omega$, we consider  a discrete graph $\Omega_D$. This means
that the vertices stay in their original positions, but each edge
is now represented by a finite set of points. In other words, each
edge $e_j$ in the set $E$ is identified not with the interval $(0,
l_j)$ on the real line, but with a finite set of numbers $a^j_0 <
a^j_1 < \ldots < a^j_{N_j}$. For simplicity, we can assume that
$a^j_0=0$, $a^j_1=1$,\ldots$a^j_{N_j}=N_j$ (Fig.1).
    \begin{figure}[h]
        \includegraphics[width=1\textwidth]{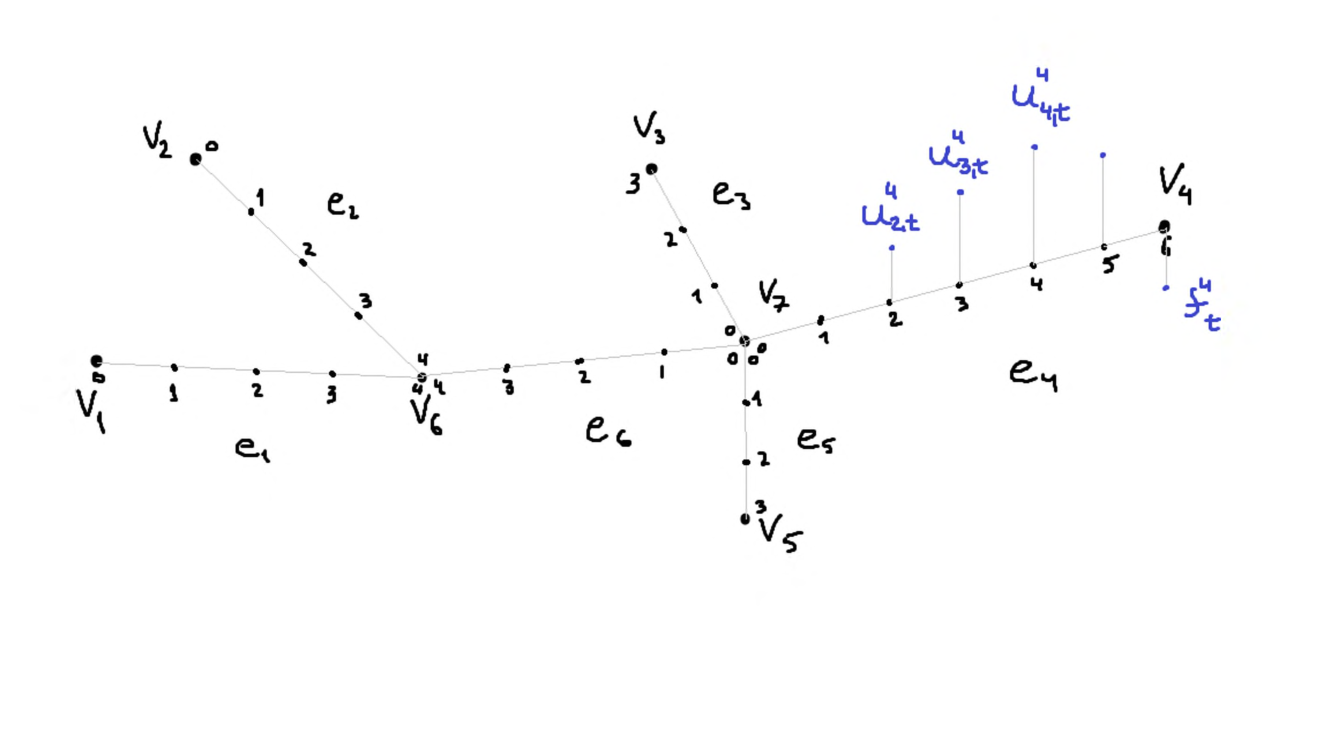}
        \caption{Discrete wave equation.}
        \label{pic2}
    \end{figure}

Additionally, we assume that the time $t$ is also discrete. In
other words,  we replace the interval $(0,T)$ with the set of
discrete points: $t_0<t_1<\dots< t_T$. For simplicity, we can
consider it as a sequence of numbers: $0<1<\ldots <T$. The space
of real square summable functions on the graph $ \Omega_D $ is
denoted by $L(\Omega_D):=  \bigoplus_{i=1}^{N}\mathbb{R}^{N_i}.$
For the function $u \in L(\Omega_D)$ we write
\begin{equation*}
u:= \left\{ u^i\right\}_{i=1}^N,\quad
u^i=(u^i_0,\ldots,u^i_{N_i}),\quad u^i_j\in \mathbb{R},\
j=0,\ldots,N_i.
\end{equation*}

In order to unify the value at the node $v_k$, we introduce the
function
\begin{equation*}
F(k)=\left[ \begin{array}l 0\\
N_k\end{array}\right.
\end{equation*}
such that
    $$
    u(v_k)=u^{j(k)}_{F(k)},\quad{\text{where } j(k)\in E(v_k)},\quad k=1,\ldots,M.
    $$

Since we are considering a discrete wave equation, the function
$u$ must depend  on two variables: coordinate and time. However,
since time is now discrete, we introduce this dependence through
an additional index on the vector $u^i$. We will express this as
follows:
    \begin{equation*}
        u(t):= \left\{ u^i(t)\right\}_{i=1}^N,\quad u^i(t)=(u^i_{0,t},\ldots,u^i_{{N_i},t}),\quad u^i_{j,t}\in \mathbb{R},\ j=0,\ldots,N_i.
    \end{equation*}
Before reformulating a direct problem on a metric graph and
deriving a similar discrete one, let us consider discrete analogs
of the kinetic and potential energy functionals:
    \begin{multline*}
        T_D(t)=\sum_{i=1}^N T^e_i(t)+\sum_{i=1}^{M}T^v_i(t)=\\
        =\sum_{i=1}^N \sum_{j=1}^{N_i-1}\frac{(u^i_{j,t}-u^i_{j,t-1} )^2}{2}  +\sum_{i=1}^{M} \frac{(u^{j(i)}_{F(i),t}-u^{j(i)}_{F(i),t-1} )^2}{2} \label{Kin_en_d},
    \end{multline*}
    \begin{equation*}
        U_D(t)=\sum_{i=1}^N U^e_i(t)= \sum_{i=1}^N \sum_{j=1}^{N_i} \frac{(u^i_{j,t}-u^i_{j-1,t} )^2}{2}. \label{Pot_en_d}
    \end{equation*}
Here, $T_i^e$, $U_i^e$, $i = 1, \ldots, N$, represent the kinetic
and potential energies associated with the interior points of the
edges $e_i$. Similarly, $T_i^v$, $i = 1, \ldots, N + 1$, denote
the kinetic energy at the vertices $v_i$. To express the action
functional in the discrete case, we need to integrate the
Lagrangian over the interval $(t_1, t_2)$. In a discrete setting,
the formula takes the following form:
\begin{equation}
\label{Discr_functional}
S[u]=\sum_{t=1}^{T}T_D(t)-
\sum_{t=0}^{T}U_D(t).
\end{equation}

Assume that the number of elements in $E(v_k)$ is $p_k$. Using the
Hamilton's principal and calculating the variations of the
functional (\ref{Discr_functional}) we come to the following
discrete system:
\begin{conclusion}
The appropriate correct formulation of the direct discrete problem
for a discrete graph would be the following system
\begin{eqnarray*}
u^i_{j,t+1}+u^i_{j,t-1}-u^i_{j+1,t}-u^i_{j-1,t}=0,\ 1\leqslant i \leqslant N,\   1\leqslant j \leqslant N_i-1, t\geqslant 0\label{Wave_eqn_d33}\\
u^i_{F(i),t}=u^j_{F(j),t},\ e_i,e_j\in E(v_k),\quad   t\geqslant 1,\label{Cont_d33}\\
u(v_k,t) = f^k(t)\quad \text{for}\,\, v_k\in\Gamma,\quad t\geqslant 1,  \label{Bound_d33}\\
u^i_{j,0}=0, \quad u^i_{j,-1}=0,\ 1\leqslant i \leqslant N,\   1\leqslant j \leqslant N_i-1\label{Init_d33}\\
-\frac{p_k}{2}u^k_{F(k),t+1}-\frac{p_k}{2}u^k_{F(k),t-1}  +
\sum_{e_i\in E(v_k)}u^{i}_{|F(i)-1|,t}=0,\label{Kirch_d_1mm3}
\end{eqnarray*}
where the last equality holds in all internal vertices.
\end{conclusion}
The derivation of the above system is given in \cite[Sections
3,4]{AChMM}.

\subsection{Application to parabolic equation.}

The idea is to substitute the system (\ref{Parab_eqn}) by the
discrete one:
\begin{equation}
\label{Heat_discr}
\begin{cases}
v_{n,t+1}-a_nv_{n+1,t}-a_{n-1}v_{n-1,t}-b_nv_{n,t}=0,\,\, t\in \mathbb{N}\cup\{0\},\,\, n\in \mathbb{N} \\
v_{n,\,0}=0,\quad n\in \mathbb{N}, \\
v_{0,\,t}=f_t,\quad 
\quad t\in \mathbb{N}\cup\{0\},
\end{cases}
\end{equation}
where $f=(f_0,f_1,\ldots)$ is a \emph{boundary control}, $a_0=1$.
Solution to (\ref{Heat_discr}) is denoted as $v^f$. We do not set
the boundary condition at the right end in (\ref{Heat_discr}) as
in (\ref{Parab_eqn}) since the speed of wave propagation in
(\ref{Heat_discr}) is finite (see Remark \ref{Rem_finite}).

For the solution of (\ref{Heat_discr}) one can derive an analog of
the representation (\ref{Jac_sol_rep}):
\begin{lemma}
\label{Heat_repr}
For the solution $v^f$ the following
representation is valid:
\begin{equation}
\label{heat_sol_rep} v^f_{n,t}=\prod_{k=0}^{n-1}
a_kf_{t-n}+\sum_{s=n}^{t-1}w_{n,s}f_{t-s-1},\quad n,t\in
\mathbb{N},
\end{equation}
where  $w_{n,s}$ satisfy the Goursat problem:
\begin{equation*}
\label{Goursat2} \left\{
\begin{array}l
w_{n,s+1}-a_nw_{n+1,s}-a_{n-1}w_{n-1,s}-b_nw_{n,s}=\delta_{s,n}a_n^2\prod_{k=0}^{n-1}a_k,\  n,s\in \mathbb{N}, \ s>n,\\
w_{n,n}-b_n\prod_{k=0}^{n-1}a_k-a_{n-1}w_{n-1,n-1}=0,\quad w_{n,n-1}=0,\quad n\in \mathbb{N},\\
w_{0,t}=0,\quad t\in \mathbb{N}_0.
\end{array}
\right.
\end{equation*}
\end{lemma}
The proof of Lemma (\ref{Heat_repr}) is given in \cite{MMSh}.

We define the \emph{response operator}: $R^T:\mathbb{R}^T\mapsto
\mathbb{R}^T$ by the rule:
\begin{equation*}
R^Tf:=(v^f_{1,1},v^f_{1,2},\ldots,v^f_{1,T} ).
\end{equation*}
Using (\ref{heat_sol_rep}) with $n=1$, we get that
\begin{equation*}\label{Reac1}
\left(R^Tf\right)_t=a_0 f_{t-1}+ \sum_{s=1}^{t-1}w_{1,s}f_{t-s-1}.
\end{equation*}
Introducing the definition \emph{response vector (convolution
kernel of the reaction operator)}\,---\,
$s=(s_0,s_1,\ldots,s_{T-1})=(a_0,w_{1,1},w_{1,2},\ldots,w_{1,T-1})$.
Setting $f=\delta=(1,0,\ldots)$, we rewrite the last equality in
the form
\begin{equation*} 
\left(R^T\delta\right)_t= v^\delta_{1,t}= s_{t-1},\quad
t=1,\ldots,T.
\end{equation*}
The dynamic IP for (\ref{Heat_discr}) is to recover a Jacobi
matrix (i.e. the sequences $\{a_1,a_2,\ldots\}$, $\{b_1,
b_2,\ldots \}$) and additional parameter $a_0$ from the response
operator $R^T$.

For the system (\ref{Heat_discr}) we introduce the \emph{control
operator} $V^T:\mathbb{R}^T\mapsto\mathbb{R}^T $
\begin{equation*}
V^Tf:=(v^f_{1,T},v^f_{2,T},\ldots v^f_{T,T}).
\end{equation*}

From (\ref{heat_sol_rep}) it follows that
\begin{equation*}\label{W_rep1}
(V^T)_n = v^f_{n,T}= \prod_{k=0}^{n-1}
a_kf_{T-n}+\sum_{s=n}^{T-1}w_{n,s}f_{T-s-1},\quad n=1,\ldots,T.
\end{equation*}

\begin{theorem}
\label{Heat_control} The operator $V^T$ is an isomorphism in the
space $\mathbb{R^T}$.
\end{theorem}

For the system (\ref{Heat_discr}) we introduce the
\emph{connecting operator} $S^T:\mathbb{R}^T\mapsto\mathbb{R}^T $
\begin{equation*}\label{ST_def}
 S^T:=(V^T)^* V^T,\quad (S^T f,g)= (V^T f,V^T g )_{\mathbb{R^T}}.
\end{equation*}

\begin{theorem}  The connecting operator $S^T$ is an isomorphism in the space $\mathbb{R}^T$. It can be represented using inverse data:
    \begin{equation*}
    \label{Heat_CT} S^T=
    \begin{pmatrix}
    s_{2T-2}&\ldots & \ldots & s_T & s_{T-1}\\
    s_{2T-3}&\ldots & \ldots & \ldots &s_{T-2}\\
    \ldots &\ldots &\ldots & \ldots & \ldots \\
    s_{T+1} &\ldots & s_4 & s_3 & s_2\\
    s_{T}& \ldots&s_3&s_2&s_1 \\
    s_{T-1}&\ldots &  s_2 & s_1 &s_0
    \end{pmatrix},\quad S^T_{ij}= s_{2T-(i+j)},
    \end{equation*}
\end{theorem}
The proofs of Theorems \ref{Heat_control} and \ref{Heat_CT} are
given in \cite{MMSh}.

Since the connecting operator for the system (\ref{Heat_discr})
has the form of a  Hankel matrix, one can use the results of
Sections \ref{Section_Jacobi}, \ref{Moment_problem} to recover the
entries of $A$.

\noindent{\bf Acknowledgments}

The authors express their deep gratitude to the anonymous reviewer
for valuable comments and suggestions.

\end{document}